\documentclass[a4paper,12pt,authoryear]{elsarticle} %,
\usepackage[left=3.0cm,right=2.0cm,top=2cm,bottom=3cm,bindingoffset=0cm]{geometry}

\usepackage{graphicx}
\usepackage{natbib}
\usepackage[english]{babel}

\usepackage{tabularx}
\usepackage{multirow}
\usepackage{multicol}
\usepackage{amssymb}
\usepackage{amstext}
\usepackage{amsmath}
\usepackage{pslatex}
\usepackage{floatrow}

\usepackage[update,prepend]{epstopdf}

\newtheorem{definition}{Definition}

\usepackage{lmodern}

\journal{arXiv}

\begin{document}

\begin{frontmatter}

\title{Lock-in range of classical PLL with impulse signals and proportionally-integrating filter}

\author{K.~D. Aleksandrov}
\author{N.V. Kuznetsov\corref{cor}\,}
\ead{nkuznetsov239@gmail.com}
\author{G.~A. Leonov}
\author{M.~V. Yuldashev}
\author{R.~V. Yuldashev}

\address{Faculty of Mathematics and Mechanics,\
Saint-Petersburg State University, Russia}
\address{Dept. of Mathematical Information Technology,\
University of Jyv\"{a}skyl\"{a}, Finland}
\address{Institute of Problems of Mechanical Engineering RAS, Russia}

\begin{abstract}
In the present work the model of PLL with impulse signals and active PI filter in the signal's phase space is described.
For the considered PLL the lock-in range is computed analytically
and obtained result are compared with numerical simulations.
\end{abstract}


\end{frontmatter}
\section{Models of classical PLL with impulse signals}
\label{sec:PLLBased}
\begin{figure}[!htbp]
\centering
\includegraphics[width=0.9\textwidth]{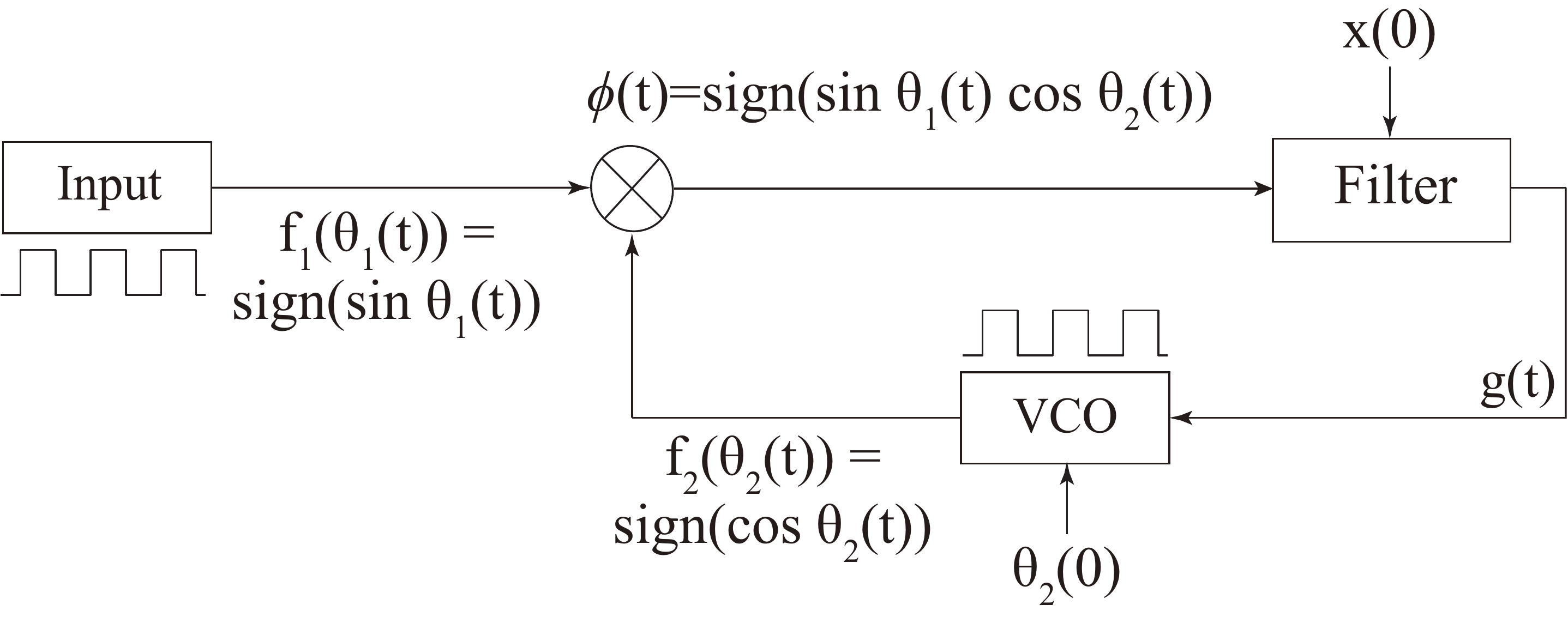}
\caption{Model of PLL with impulse signals in the signals space.}
\label{ris:PLLSigSpace}
\end{figure}

Consider a physical model of classical PLL in the signals space (see Fig.~\ref{ris:PLLSigSpace}).
% Models of the PLL in the signals space are difficult for the study \citep{KudrewiczW-2007} since the equations, which describe these models, are nonautonomous. By contrast, equations for the models in the signal's phase space are autonomous \citep{Viterbi-1966, Shahgildyan-1966, Gardner-1966}, what simplifies their study.

% From the numerical point of view, advantage of models in the signal's phase space is the
% nonexistence of high-frequency components, thus simulation in the signal's phase space allows one to consider
% slow varying frequency only.
% By contrast, the simulation of PLL in
% the signals space is complicated
% since we have to observe simultaneously both high-frequency
% (fast changing of phases) and low-frequency (relatively slow changing
% of frequencies) oscillations.
This model contains the following blocks: a reference oscillator (Input), a voltage-controlled oscillator (VCO), a filter (Filter), and an analog multiplier as a phase detector (PD).
The signals $\operatorname{sign}\left(\sin \theta_1(t)\right)$ and $\operatorname{sign}\left(\cos \theta_2(t)\right)$ of the Input and the VCO (here $\theta_2(0)$ is the initial phase of VCO) enter the multiplier block. The resulting impulse signal
$\phi(t) = \operatorname{sign}\left(\sin \theta_1(t) \cos \theta_2(t)\right)$ is filtered by low-pass filter Filter (here $x(0)$ is an initial state of Filter). The filtered signal $g(t)$ is used as a control signal for VCO.

The equations describing the model of PLL-based circuits in the signals space are difficult for the study, since that equations are nonautonomous (see, e.g., \citep{KudrewiczW-2007}). By contrast, the equations of model in the signal's phase space are autonomous
\citep{Gardner-1966,ShahgildyanL-1966,Viterbi-1966}, what simplifies the study of PLL-based circuits.
The application of averaging methods \citep{MitropolskyB-1961,Samoilenko-2004-averiging} allows one to reduce the model of PLL-based circuits in the signals space to the model in
the signal's phase space
(see, e.g., \citep{LeonovKYY-2012-TCASII,LeonovK-2014-book,LeonovKYY-2015-SIGPRO,KuznetsovLSYY-2015-PD,KuznetsovKLNYY-2015-ISCAS,BestKKLYY-2015-ACC}.
\begin{figure}[!htbp]
\centering
\includegraphics[width=0.9\textwidth]{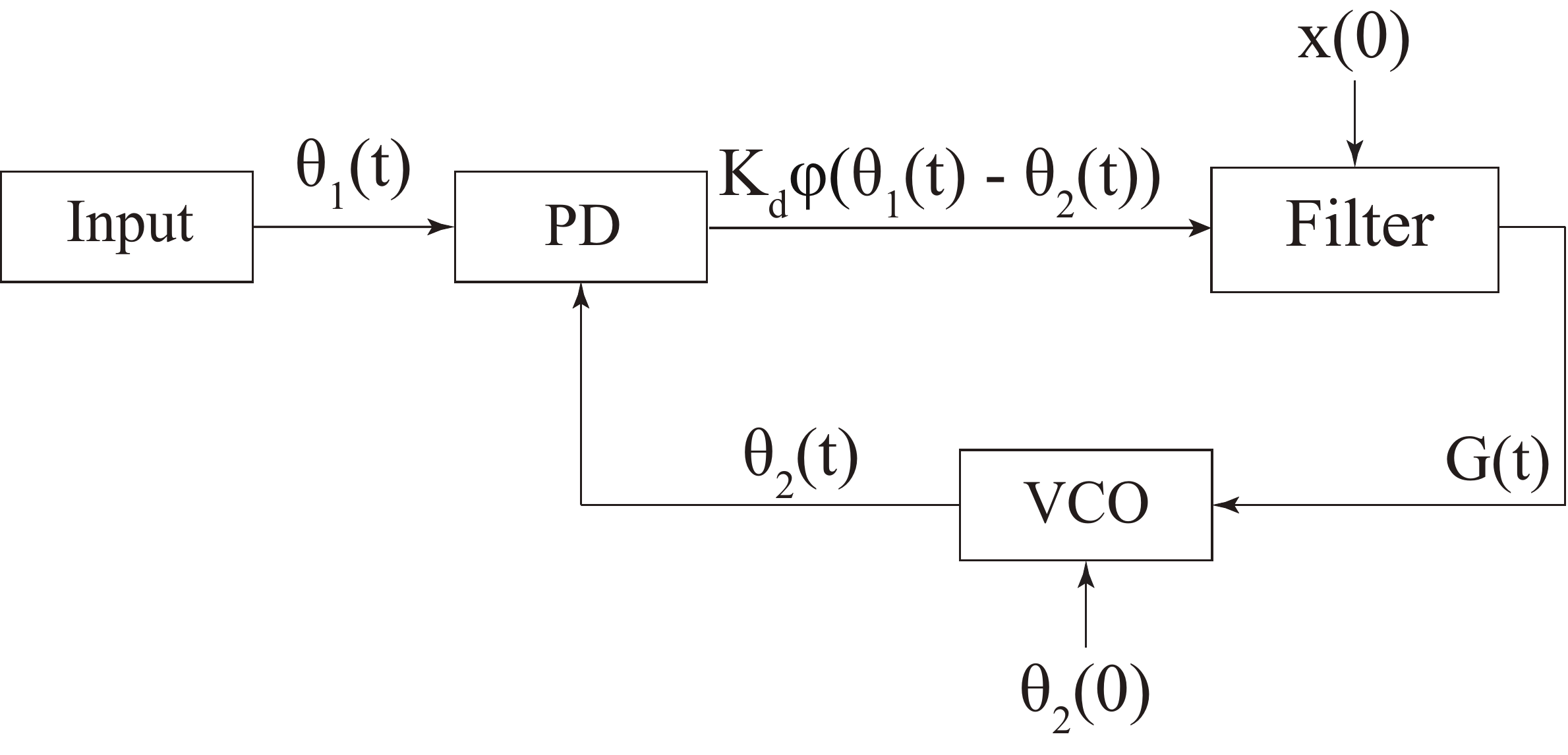}
\caption{Model of the classical PLL in the signal's phase space.}
\label{ris:PLL_PF_space}
\end{figure}

The main difference between the physical model (Fig.~\ref{ris:PLLSigSpace}) and the simplified mathematical model in the signal's phase space (Fig.~\ref{ris:PLL_PF_space}) is the absence of high-frequency component of the phase detector output.
% The high-frequency component
% \begin{equation}
% \psi(t) = \phi(t) - K_d \varphi(\theta_\Delta(t))
% \label{rel:PDSignals}
% \end{equation}
% is suppressed by the Filter.
The output of the phase detector in the signal's phase space is called a phase detector characteristic and has the form $$ K_d \varphi(\theta_1(t) - \theta_2(t)).$$ The maximum absolute value of PD output $K_d > 0$ is called a phase detector gain (see, e.g., \citep{Best-2007, Goldman-2007-book}). The periodic function $\varphi(\theta_\Delta(t))$ depends on difference $\theta_1(t) - \theta_2(t)$ (which is called a phase error and denoted by $\theta_\Delta(t)$).
The PD characteristic depends on the design of PLL-based circuit and the signal waveforms $f_1(\theta_1)$ of Input and $f_2(\theta_2)$ of VCO.
For PLL with impulse signals the PD characteristic is as follows (see, e.g., \citep{Viterbi-1966,Gardner-1966,LeonovKYY-2012-TCASII}):
\begin{align}
& K_d = 1; \nonumber \\
&\varphi(\theta_\Delta(t))=
\begin{cases}
\frac{2}{\pi}\theta_\Delta(t), &\text{if $-\frac{\pi}{2} \leq \theta_\Delta(t) \leq \frac{\pi}{2}$,}\\
-\frac{2}{\pi}\theta_\Delta(t) + 2, &\text{if $\frac{\pi}{2} \leq \theta_\Delta(t) \leq \frac{3\pi}{2}$}.
\end{cases}
\label{eq:PDCharTriangle}
\end{align}

Let us describe a model of classical PLL with impulse signals in the signal's phase
space (see Fig.~\ref{ris:PLL_PF_space}).
A reference oscillator and a voltage-controlled oscillator generate the phases
$\theta_1(t)$ and $\theta_2(t)$, respectively.
The frequency of reference signal usually assumed to be constant:
\begin{equation}
\dot{\theta}_1(t) = \omega_1.
\label{eq:Input}
\end{equation}

The phases $\theta_1(t)$ and $\theta_2(t)$ enter the inputs of the phase detector.
The output of phase detector is processed by Filter. Further we consider the active PI filter (see, e.g., \citep{Baker-2011-book}) with transfer function
$W(s) = \frac{1 + \tau_2 s}{\tau_1 s}$, $\tau_1~>~0$, $\tau_2~>~0$. The considered filter can be described as
% \begin{equation}\dot{G}(t) = \frac{\tau_2 K_d}{\tau_1}\cos(\theta_\Delta(t))\dot{\theta}_\Delta(t) + \frac{K_d}{\tau_1} \sin(\theta_\Delta(t))
% \label{eq:FilterOrig}
% \end{equation}
% \begin{equation}
% \begin{cases}
% \tau_1\dot{x}(t) = K_d\varphi(\theta_\Delta(t)), \\
% G(t) = \frac{\tau_2}{\tau_1}K_d\varphi(\theta_\Delta(t)) + x(t).
% \label{eq:FilterOrig0}
% \end{cases}
% \end{equation}
% System (\ref{eq:FilterOrig0}) can be rewritten in the equivalent form:
\begin{equation}
\begin{cases}
\dot{x}(t) = K_d\varphi(\theta_\Delta(t)), \\
G(t) = \frac{1}{\tau_1}x(t) + \frac{\tau_2}{\tau_1}K_d\varphi(\theta_\Delta(t)),
\label{eq:FilterOrig}
\end{cases}
\end{equation}
where $x(t)$ is the filter state.

% The Filter with transfer function is described by the equation
% \begin{equation}
% \dot{G}(t) = \frac{\tau_2 K_d}{\tau_1}\cos(\theta_\Delta(t))\dot{\theta}_\Delta(t) + \frac{K_d}{\tau_1} \sin(\theta_\Delta(t))
% \label{eq:FilterOrig}
% \end{equation}
% where $x(t)$ is filter state.
% By linear change $x(t) \rightarrow \frac{K_d}{\tau_1}\left(x(t) + \tau_2 \sin(\theta_\Delta(t))\right)$ (\ref{eq:FilterOrig}) takes the form:
% \begin{equation}
% \begin{cases}
% \dot{x}(t) = \sin(\theta_\Delta(t)), \\
% \frac{K_d}{\tau_1}\left(x(t) + \tau_2 \sin(\theta_\Delta(t))\right) = G(t).
% \end{cases}
% \label{eq:Filter}
% \end{equation}
The output of Filter $G(t)$ is used as a control signal for VCO:
\begin{equation}
\dot{\theta}_2(t) = \omega_2^{\rm free} + K_v G(t),
\label{eq:VCO}
\end{equation}
where $\omega_2^{\rm free}$ is the VCO free-running frequency and $K_v > 0$ is the VCO gain.

Relations (\ref{eq:Input}), (\ref{eq:FilterOrig}), and (\ref{eq:VCO}) result in autonomous system of differential equations
\begin{equation}
\begin{cases}
\dot{x} = K_d \varphi(\theta_\Delta), \\
\dot{\theta}_\Delta = \omega_1 - \omega_2^{\rm free} - \frac{K_v}{\tau_1}\left(x + \tau_2 K_d \varphi(\theta_\Delta)\right).
\end{cases}
\label{sys:PLLSys_beforeChange}
\end{equation}
Denote the difference of the reference frequency and the VCO free-running frequency $\omega_1 - \omega_2^{\rm free}$ by $\omega_\Delta^{\rm free}$.
By the linear transformation $x \rightarrow K_d x$ we have
\begin{equation}
\begin{cases}
\dot{x} = \varphi(\theta_\Delta), \\
\dot{\theta}_\Delta = \omega_\Delta^{\rm free} - \frac{K_0}{\tau_1}\left(x + \tau_2 \varphi(\theta_\Delta)\right),
\end{cases}
\label{sys:PLLSys}
\end{equation}
where $K_0 = K_v K_d$ is the loop gain. Here (\ref{sys:PLLSys}) describes the model of PLL with the impulse signals and active PI filter in the signal's phase space.

By the transformation $$\left(\omega_\Delta^{\rm free}, x, \theta_\Delta\right) \rightarrow \left(-\omega_\Delta^{\rm free}, -x, -\theta_\Delta\right),$$ (\ref{sys:PLLSys}) with odd PD characteristic (\ref{eq:PDCharTriangle}) is not changed.
This property allows one to use the concept of frequency deviation $$\left|\omega_\Delta^{\rm free}\right| = \left|\omega_1 - \omega_2^{\rm free}\right|$$ and consider (\ref{sys:PLLSys}) with $\omega_\Delta^{\rm free} > 0$ only.

The PLL state for which the VCO frequency is adjusted to the reference frequency of Input is called a locked state.
The locked states of the PLL correspond to the locally asymptotically stable equilibria of (\ref{sys:PLLSys}), which can be found from the relations
\begin{equation}
\begin{cases}
\varphi(\theta_{eq}) = 0, \\
\omega_\Delta^{\rm free} - \frac{K_0}{\tau_1} x_{eq} = 0.
\end{cases}
\nonumber
\label{sys:PLLSysEq}
\end{equation}

Since (\ref{sys:PLLSys}) is $2\pi$-periodic in $\theta_\Delta$, we can consider (\ref{sys:PLLSys}) in a $2\pi$-interval of $\theta_\Delta$, $\theta_\Delta \in \left(-\pi, \pi\right]$.
In interval $\theta_\Delta \in \left(-\pi, \pi\right]$ there exist two equilibria: $$\left(\theta_{eq}^s, x_{eq}(\omega_\Delta^{\rm free})\right) = (0, \frac{\omega_\Delta^{\rm free} \tau_1}{K_0}) \text{ and } \left(\theta_{eq}^u, x_{eq}(\omega_\Delta^{\rm free})\right) = (\pi, \frac{\omega_\Delta^{\rm free} \tau_1}{K_0}).$$
As is shown below (see \ref{sec:AppLockIn}) the equilibria $$\left(\theta_{eq}^s + 2\pi k, x_{eq}(\omega_\Delta^{\rm free})\right) = \left(2\pi k, \frac{\omega_\Delta^{\rm free} \tau_1}{K_0}\right)$$ are locally asymptotically stable. Hence, the locked states of (\ref{sys:PLLSys}) are given by equilibria $\left(\theta_{eq}^s, x_{eq}(\omega_\Delta^{\rm free})\right)$. The remaining equilibria $$\left(\theta_{eq}^u + 2\pi k, x_{eq}(\omega_\Delta^{\rm free})\right) = \left(\pi + 2\pi k, \frac{\omega_\Delta^{\rm free} \tau_1}{K_0}\right)$$ are saddle equilibria (see \ref{sec:AppLockIn}).
\section{The lock-in range}
\label{sec:LockInDef}
\begin{figure}[!htbp]
\includegraphics[width=\linewidth]{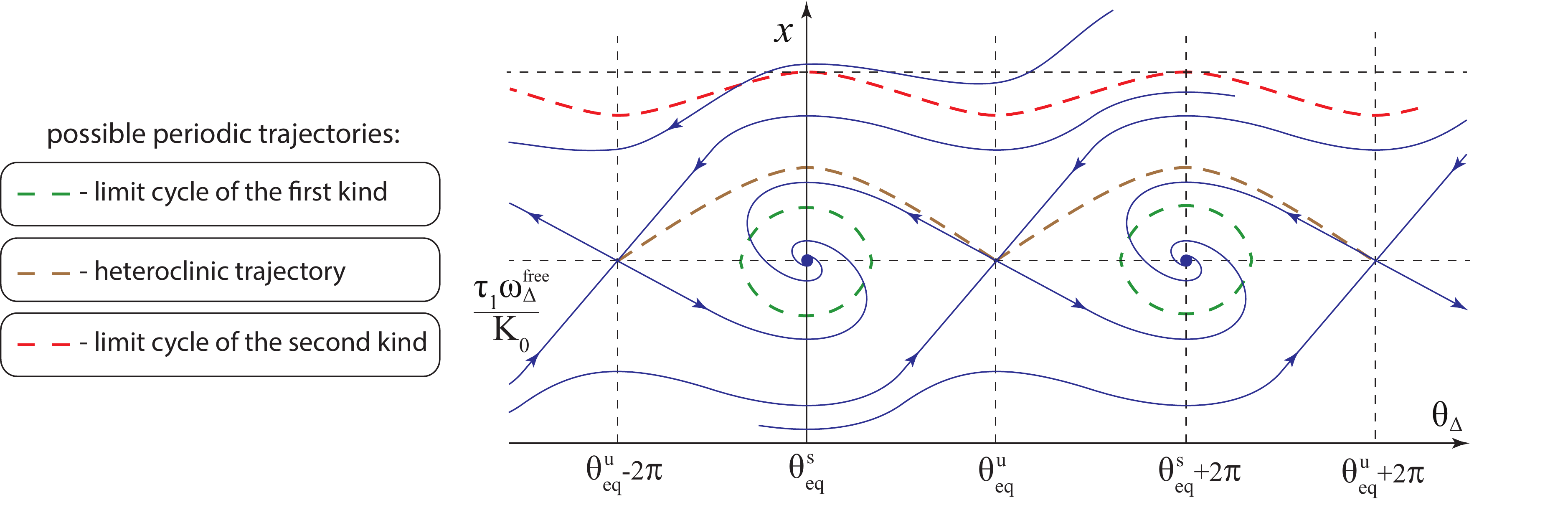}
\caption{Possible periodic trajectories on the phase plane of (\ref{sys:PLLSys}).}
\label{ris:CyclesPhasePlane}
\end{figure}

The model of classical PLL with impulse signals and active PI filter in the signal's phase space is globally asymptotically stable (see, e.g., \citep{Gubar-1961,LeonovA-DAN-2015}). The PLL achieves locked state for any initial VCO phase $\theta_2(0)$ and filter state $x(0)$. So, there exist no limit cycles of the first kind, heteroclinic trajectories, and limit cycles of the second kind on the phase plane of (\ref{sys:PLLSys}) (see Fig.~\ref{ris:CyclesPhasePlane}).

However, the phase error $\theta_\Delta$ may significantly increase during the acquisition process. In order to consider the property of the model to synchronize without undesired growth of the phase error $\theta_\Delta$, a lock-in range concept was introduced in \citep{Gardner-1966}:
``{\it{If, for some reason,  the frequency difference between input and VCO is less than the loop bandwidth, the loop will lock up almost instantaneously without slipping cycles. The maximum frequency difference for which this fast acquisition is possible is called the lock-in frequency}}''.
The lock-in range concept is widely used in engineering literature on the PLL-based circuits study (see, e.g., \citep{Stensby-1997,KiharaOE-2002,Kroupa-2003,Gardner-2005-book,Best-2007}).
It is said that a cycle slipping occurs if (see, e.g., \citep{AscheidM-1982,ErshovaL-1983,SmirnovaPU-2014})
\begin{equation}
\displaystyle \limsup_{t \rightarrow +\infty} \left|\theta_\Delta(0) - \theta_\Delta(t) \right| \geq 2\pi.
\nonumber
\end{equation}
\begin{figure}[!htbp]
\includegraphics[width=\linewidth]{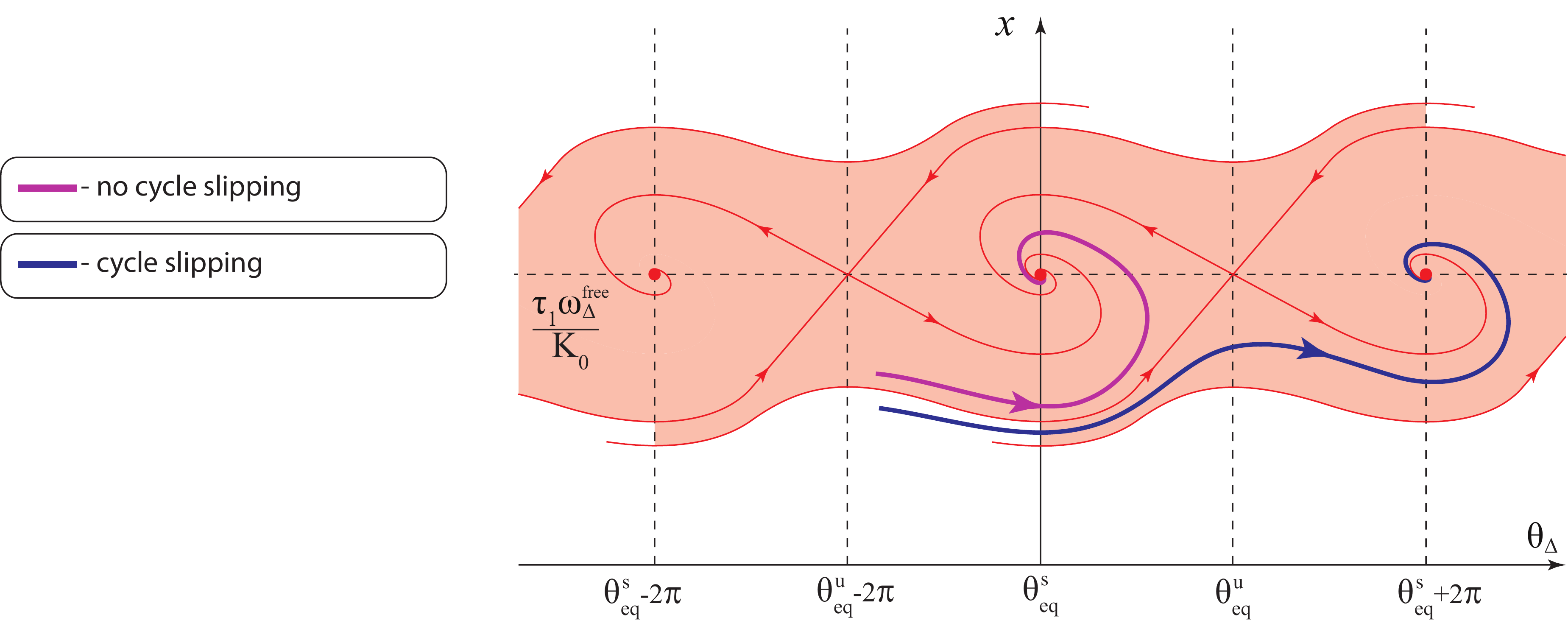}
\caption{The lock-in domain and cycle slipping.}
\label{ris:CycleSlipping}
\end{figure}
However, in general, even for zero frequency deviation ($\omega_\Delta^{\text{free}}=0$)
and a sufficiently large initial state of filter ($x(0)$),
cycle slipping may take place, thus in 1979 Gardner wrote: \textit{``There is no natural way to define exactly any unique lock-in frequency''} and \textit{``despite its vague reality, lock-in range is a useful concept''} \citep{Gardner-1979-book}.

To overcome the stated problem, in \citep{KuznetsovLYY-2015-IFAC-Ranges,LeonovKYY-2015-TCAS} the rigorous mathematical definition of a lock-in range is suggested:

% The pull-in range from the definition of the lock-in range is \textit{the largest interval $[0, \omega_p)$ of frequency deviations $|\omega_\Delta^{\rm free}|$ such that the mathematical model of the loop in the signal's phase space is globally asymptotically stable} \citep{LeonovKYY-TCASI-2015}.
\begin{definition} \citep{KuznetsovLYY-2015-IFAC-Ranges,LeonovKYY-2015-TCAS}
\it The lock-in range of model (\ref{sys:PLLSys}) is a range $\left[0, \omega_l\right)$ such that for each frequency deviation $\left|\omega_\Delta^{\rm free}\right| \in \left[0, \omega_l\right)$ the model (\ref{sys:PLLSys}) is globally asymptotically stable and the following domain $$D_{\rm lock-in}\left((-\omega_l, \omega_l)\right) = \underset{\left|\omega_\Delta^{\rm free}\right|<\omega_l}{\bigcap} D_{\rm lock-in}(\omega_\Delta^{\rm free})$$ contains all corresponding equilibria $\left(\theta^s_{eq}, x_{eq}(\omega_\Delta^{\rm free})\right).$
\label{def:LockIn}
\end{definition}
For model (\ref{sys:PLLSys}) each lock-in domain from intersection $\underset{\left|\omega_\Delta^{\rm free}\right|<\omega_l}{\bigcap} D_{\rm lock-in}(\omega_\Delta^{\rm free})$ is bounded by the separatrices of saddle equilibria $\left(\theta^u_{eq}, x_{eq}(\omega_\Delta^{\rm free})\right)$ and vertical lines $\theta_\Delta = \theta^s_{eq} \pm 2\pi$.
Thus, the behavior of separatrices on the phase plane is the key to the lock-in range study (see Fig.~\ref{ris:LockinPhasePlane}).
% Two corresponding separatrices of each saddle equilibrium of (\ref{sys:PLLSys}) are the only phase trajectories, which do not tend to the stable equilibria. Moreover, those separatrices separate phase trajectories, which tend to the neighboring stable equilibria (see Fig.~\ref{ris:CyclesPhasePlane}). Thus, to determine if cycle slips occur, it is necessary to study behavior of separatrices on the phase plane.
\section{Phase plane analysis for the lock-in range estimation}
\label{sec:PhasePlaneAnalysis}

\begin{figure}[!htbp]
\includegraphics[width=0.7\linewidth]{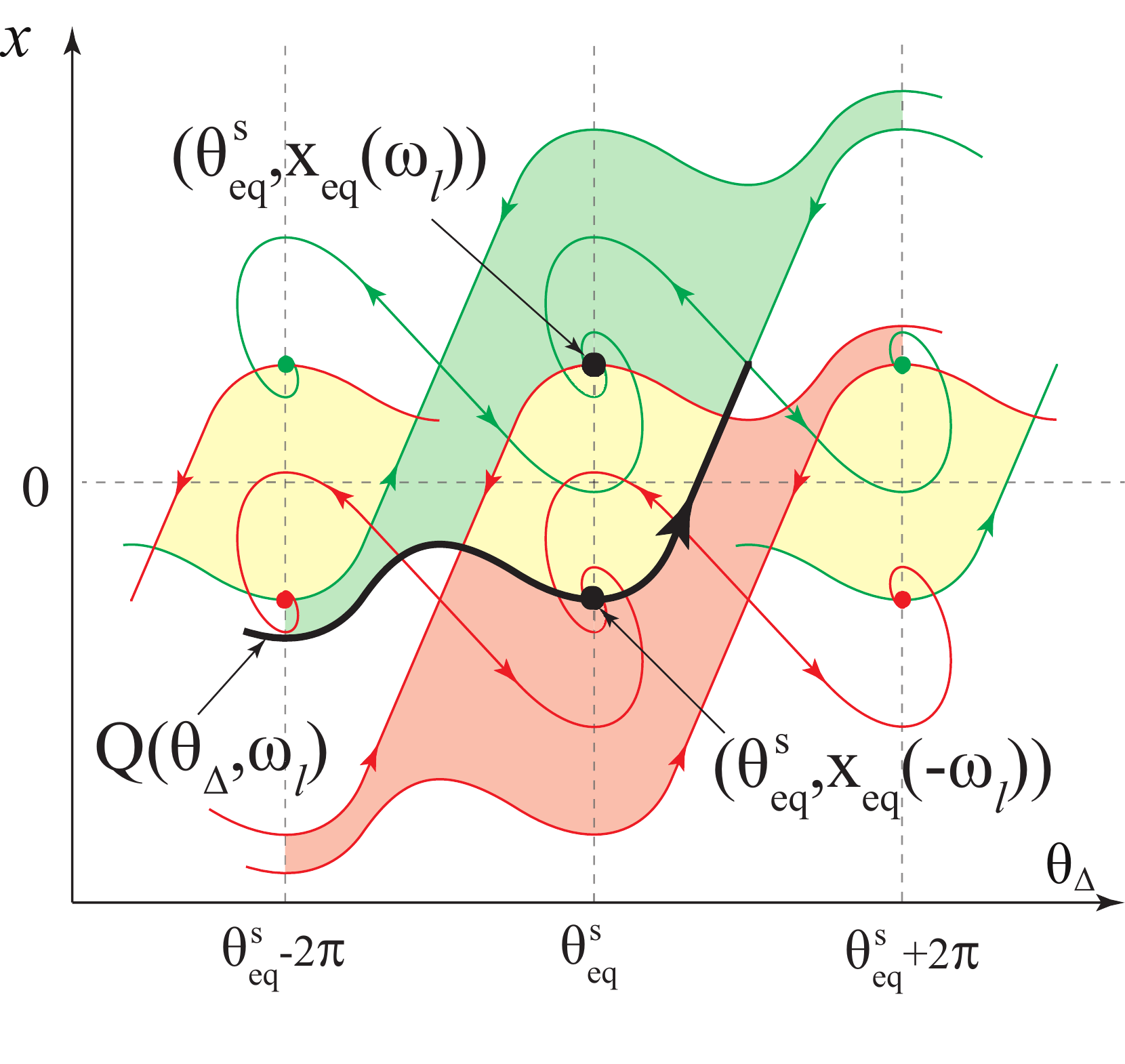}
\caption{The lock-in domain of (\ref{sys:PLLSys}) for $\left|\omega_\Delta^{\rm free}\right| = \omega_l$.}
\label{ris:LockinPhasePlane}
\end{figure}

Consider an approach to the lock-in range computation of (\ref{sys:PLLSys}), based on the phase plane analysis.
% To study the phase plane of the (\ref{sys:PLLSys}) we can consider the behavior of the phase plane trajectories on interval $\theta_\Delta \in \left[-\pi, \pi\right]$ only, since the (\ref{sys:PLLSys}) is periodic on variable $\theta_\Delta$.
To compute the lock-in range of (\ref{sys:PLLSys}) we need to consider the behavior of the lower separatrix $Q(\theta_\Delta, \omega_\Delta^{\rm free})$, which tends to the saddle point $\left(\theta_{eq}^u, x_{eq}(\omega_\Delta^{\rm free})\right) = \left(\pi, \frac{\omega_\Delta^{\rm free} \tau_1}{K_0}\right)$ as $t \rightarrow +\infty$ (by the symmetry of the lower and the upper half-planes, the consideration of the upper separatrix is also possible).
The parameter $\omega_\Delta^{\rm free}$ shifts the phase plane vertically. To check this, we use a linear transformation $x \rightarrow x + \frac{\omega_\Delta^{\rm free} \tau_1}{K_0}$. Thus, to compute the lock-in range of (\ref{sys:PLLSys}), we need to find $\omega_\Delta^{\rm free} = \omega_l$ (where $\omega_l$ is called a lock-in frequency) such that (see Fig.~\ref{ris:LockinPhasePlane})
\begin{equation}
x_{eq}(-\omega_l) = Q(\theta^s_{eq}, \omega_l).
\label{eq:LockinRelation}
\end{equation}
By (\ref{eq:LockinRelation}), we obtain an exact formula for the lock-in frequency $\omega_l$:
\begin{align}
&-\frac{\omega_l}{K_0 / \tau_1} = \frac{\omega_l}{K_0 / \tau_1} + Q(\theta^s_{eq}, 0). \nonumber \\
&\omega_l = -\frac{K_0 Q(\theta^s_{eq}, 0)}{2 \tau_1},
\label{eq:LockinRelationNew}
\end{align}

\begin{figure}[!htbp]
\centering
\includegraphics[width=0.9\textwidth]{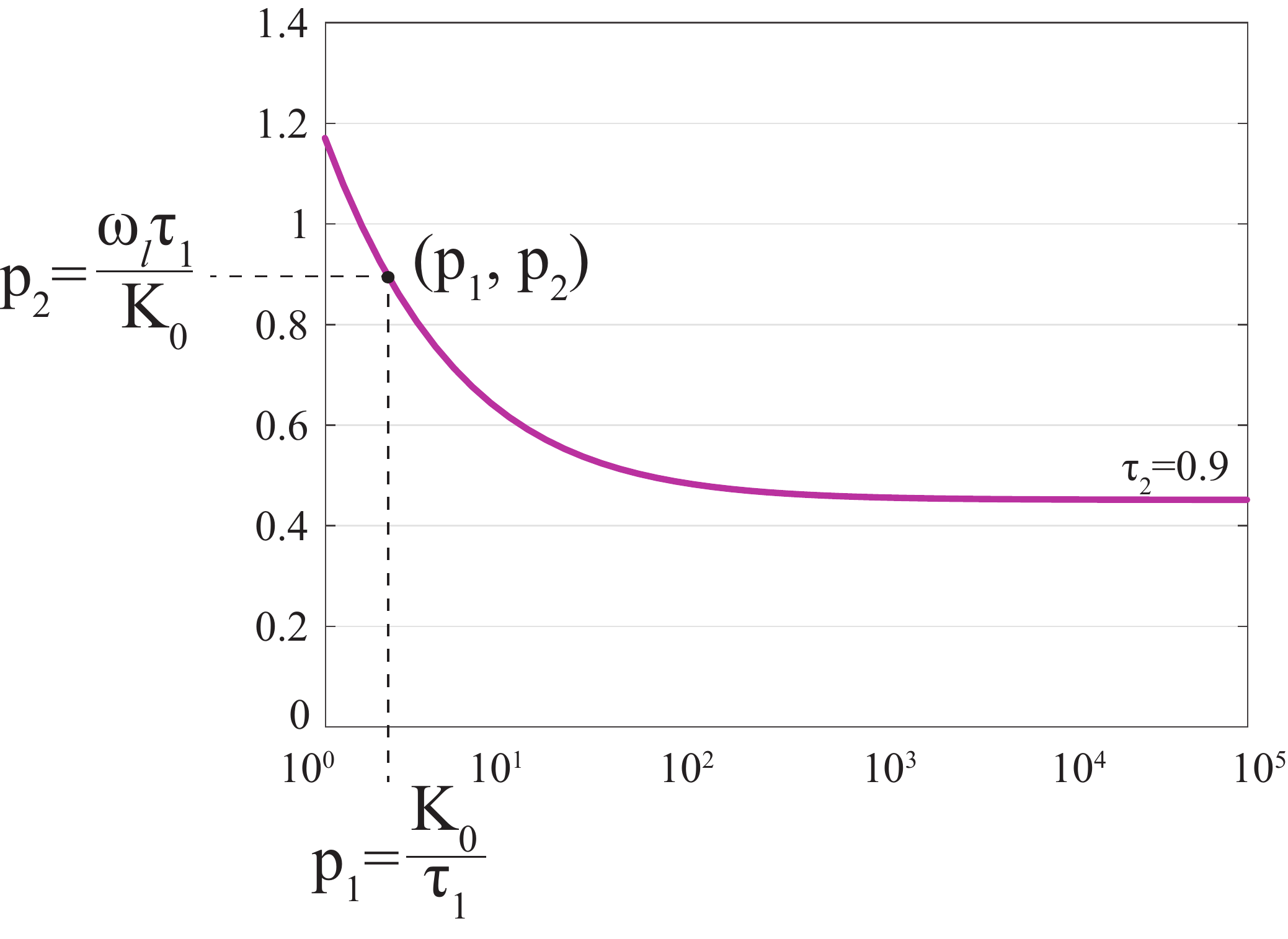}
\caption{Diagram for the lock-in frequency $\omega_l$ calculation.}
\label{ris:DiagramTriangleExample}
\end{figure}

Numerical simulations are used to compute the lock-in range of (\ref{sys:PLLSys}) applying (\ref{eq:LockinRelationNew}). The separatrix $Q(\theta_\Delta, 0)$ is numerically integrated and the corresponding $\omega_l$ is approximated. The obtained numerical results can be illustrated by special diagram (see Fig.~\ref{ris:DiagramTriangleExample}).
Note that (\ref{sys:PLLSys}) depends on the value of two coefficients $\frac{K_0}{\tau_1}$ and $\tau_2$. In Fig.~\ref{ris:DiagramTriangleExample}, choosing X-axis as $\frac{K_0}{\tau_1}$, we can plot a single curve for every fixed value of $\tau_2$. The results of numerical simulations show that for sufficiently large $\frac{K_0}{\tau_1}$, the value of $\omega_l$ grows almost proportionally to $\frac{K_0}{\tau_1}$. Hence, $\frac{\omega_l \tau_1}{K_0}$ is almost constant for sufficiently large $\frac{K_0}{\tau_1}$ and in Fig.~\ref{ris:DiagramTriangleExample} the Y-axis can be chosen as $\frac{\omega_\Delta^{\rm free} \tau_1}{K_0}$.

To obtain the lock-in frequency $\omega_l$ for fixed $\tau_1$, $\tau_2$, and $K_0$ using Fig.~\ref{ris:DiagramTriangleExample}, we consider the curve corresponding to the chosen $\tau_2$. Next, for X-value equal $\frac{K_0}{\tau_1}$ we get the Y-value of the curve. Finally, we multiply the Y-value by $\frac{K_0}{\tau_1}$.

Consider an analytical approach to the exact lock-in range computation.
Main stages of computation are presented in Subsection \ref{subsec:AnalyticalApproach}.

\subsection{Analytical approach to the lock-in range computation}
\label{subsec:AnalyticalApproach}
Consider a system
\begin{equation}
\begin{cases}
\dot{\theta}_\Delta(t) = y(t),\\
\dot{y}(t) = - \frac{K_0\tau_2}{\tau_1}\dot{\varphi}(\theta_\Delta(t))y(t) - \frac{K_0}{\tau_1}\varphi(\theta_\Delta(t)),
\label{sys:PLLSysEquiv}
\end{cases}
\end{equation}
where $y(t) = \omega_\Delta^{\rm free} - \frac{K_0}{\tau_1}\left(x(t) + \tau_2 \varphi(\theta_\Delta(t))\right)$. Relations (\ref{sys:PLLSysEquiv}) are equivalent to (\ref{sys:PLLSys}) and allow one to exclude $\omega_\Delta^{\rm free}$ from the computation. Note that equilibria $\left(\theta_{eq}, y_{eq}\right)$ of (\ref{sys:PLLSysEquiv}) and the corresponding equilibria $\left(\theta_{eq}, x_{eq}\right)$ of (\ref{sys:PLLSys}) are of the same type and related as
\begin{equation}
\left(\theta_{eq}, y_{eq}\right) = \left(\theta_{eq}, \omega_\Delta^{\rm free} - K_0 b x_{eq}\right).
\nonumber
\end{equation}

\begin{figure}[!htbp]
\centering
\includegraphics[width=\textwidth]{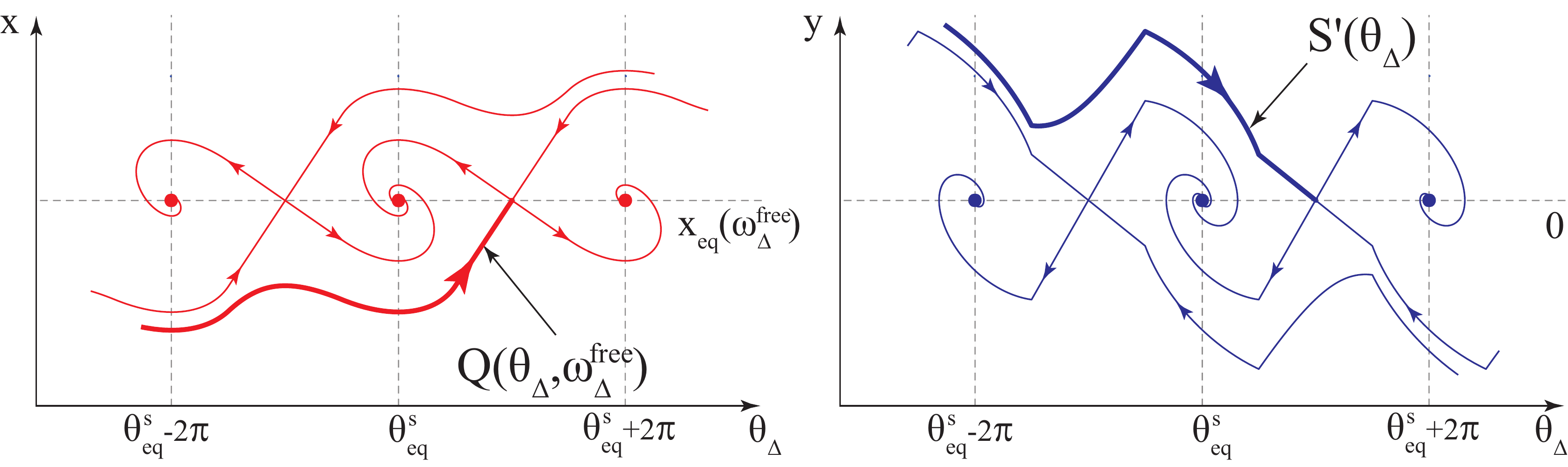}
\caption{Phase plane portraits of (\ref{sys:PLLSys}) and (\ref{sys:PLLSysEquiv}).}
\label{ris:SeparatricesEquivalent}
\end{figure}

The separatrix $Q(\theta_\Delta, \omega_\Delta^{\rm free})$ from (\ref{eq:LockinRelationNew}) corresponds to the upper separatrix $S^\prime(\theta_\Delta)$ of the phase plane of (\ref{sys:PLLSysEquiv}) (see Fig.~\ref{ris:SeparatricesEquivalent}) and the following relation
\begin{equation}
Q(\theta^s_{eq}, \omega_\Delta^{\rm free}) = \frac{\tau_1}{K_0} \left(\omega_\Delta^{\rm free} - S^\prime(\theta^s_{eq})\right) \nonumber
\end{equation}
is valid.
\\
Relation (\ref{eq:LockinRelationNew}) takes the form
\begin{equation}
\omega_l = \frac{1}{2}S^\prime(\theta^s_{eq}).
\label{eq:LockInEquivSys}
\end{equation}
% This also is the reason for consideration of (\ref{sys:PLLSysEquiv}).

The computation of the separatrix $S^\prime(\theta_\Delta)$ is in two steps. Step $1$: we integrate the separatrix $S^\prime(\theta_\Delta)$ in the interval $\left(\frac{\pi}{2}, \pi\right)$ (in which the function $\varphi(\theta_\Delta)$ is continuously differentiable) and compute $S^\prime(\frac{\pi}{2})$. For this purpose, we need to find the eigenvector that corresponds to separatrix $S^\prime(\theta_\Delta)$ on the considered interval.
Step 2: we find a general solution of (\ref{sys:PLLSysEquiv}) on the interval $\left(-\frac{\pi}{2}, \frac{\pi}{2}\right)$. Here there exist three cases depending on the type stable equilibrium $\left(\theta^s_{eq}, 0\right)$: a stable focus, stable node, and stable degenerated node. For every case described above we perform separate computations.
Using the computed $S^\prime(\frac{\pi}{2})$ as the initial data of the Cauchy problem, it is possible to obtain an exact expression for $S^\prime(\theta^s_{eq})$.

\begin{figure}[!htbp]
\centering
\includegraphics[width=0.9\textwidth]{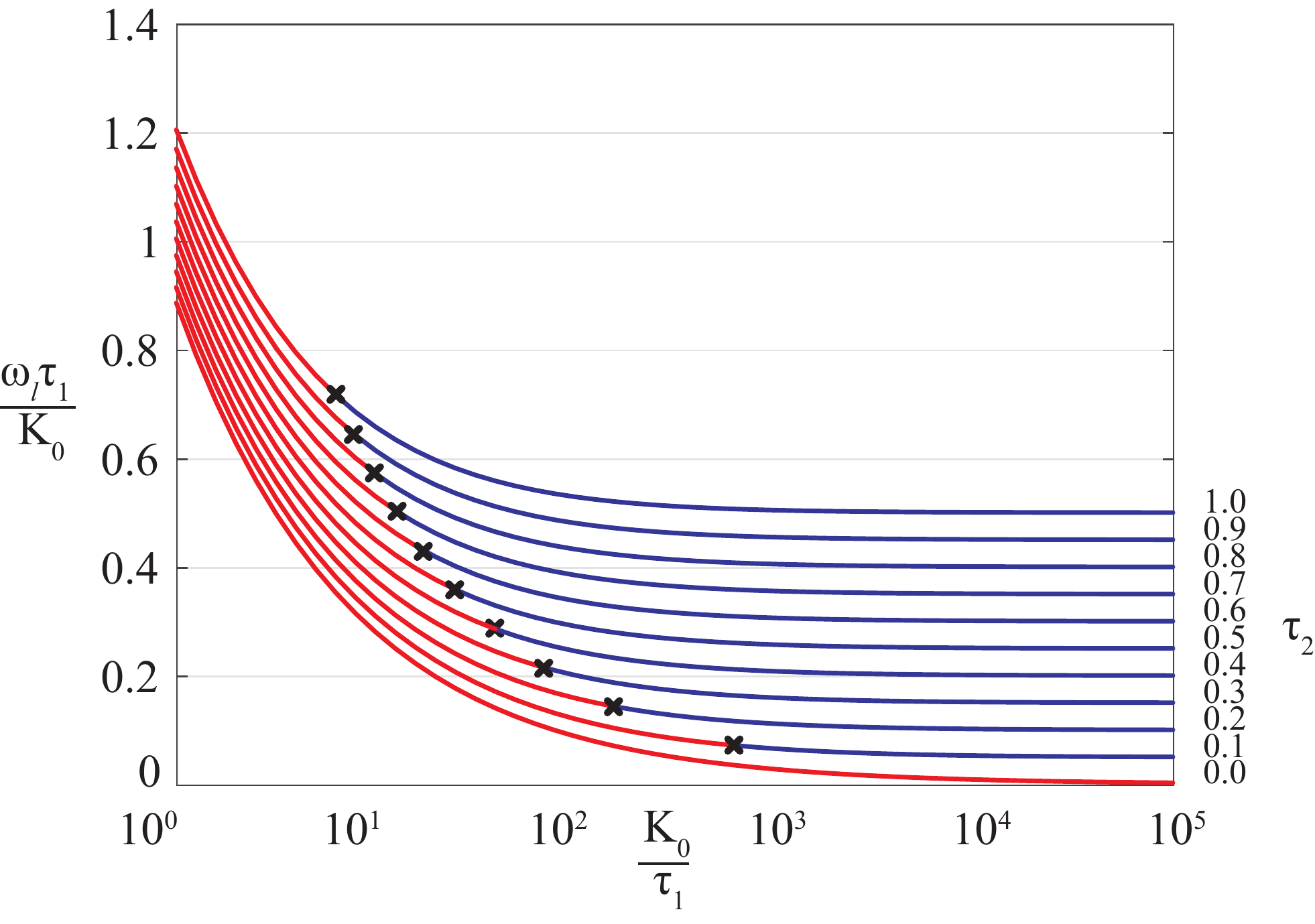}
\caption{Diagram for the lock-in frequency $\omega_l$ calculation.}
\label{ris:DiagramTriangle}
\end{figure}

The obtained analytical results are illustrated in Fig.~\ref{ris:DiagramTriangle}. The red line in Fig.~\ref{ris:DiagramTriangle} is used for the case of stable focus, and the green line for the case of stable node. The crosses are used for the case of stable degenerated node.

\begin{figure}[!htbp]
\centering
\includegraphics[width=\textwidth]{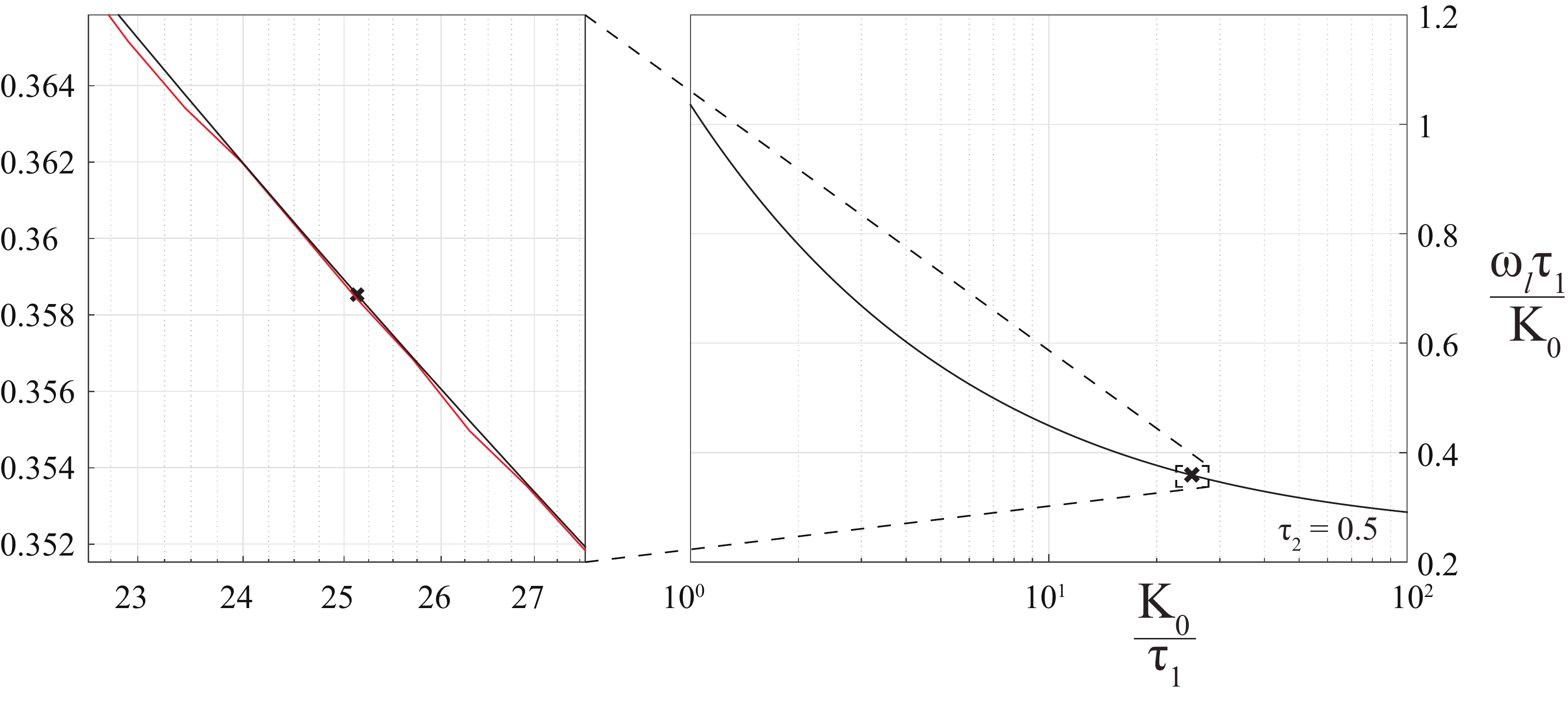}
\caption{Comparison of analytical and numerical results on the lock-in computation.}
\label{ris:DiagramTriangleComparison}
\end{figure}

The formulae for three possible cases are given below (redefinitions $a = \frac{\tau_2}{\tau_1}$, $b = \frac{1}{\tau_1}$ are used to reduce the analytical formulae): \\
\textbf{A.} $(a K_0)^2 - 2 b K_0\pi > 0$ that corresponds to a stable node:
\begin{align}
&\omega_l = \displaystyle \frac{1}{\pi} c_1 \sqrt{(a K_0)^2 - 2 b K_0\pi} \hskip0.1cm \left( -\frac{c_2}{c_1}\right)^{\left(\displaystyle \frac{1}{2} - \frac{a K_0}{\displaystyle 2 \sqrt{(a K_0)^2 - 2 b K_0\pi}}\right)}, \label{formula:LockInNode} \\
&\text{where } c_1  = \frac{\pi}{4} \left(\displaystyle \frac{\sqrt{(a K_0)^2 + 2 b K_0\pi}}{\sqrt{(a K_0)^2 - 2 b K_0\pi}} + 1\right), c_2 = \frac{\pi}{4} \left(1 - \displaystyle \frac{\sqrt{(a K_0)^2 + 2 b K_0\pi}}{\sqrt{(a K_0)^2 - 2 b K_0\pi}}\right) \nonumber.
\end{align}
\textbf{B.} $(a K_0)^2 - 2 b K_0\pi = 0$ that corresponds to a stable degenerated node:
\begin{align}
&\omega_l = \frac{1}{2}c_2 \hskip0.1cm e^{\left(\displaystyle \frac{a K_0}{2 c_2}\right)}, \text{where } c_2 = \displaystyle \frac{\sqrt{(a K_0)^2 + 2 b K_0\pi}}{2}. \label{formula:LockInDegNode}
\end{align}
\textbf{C.} $(a K_0)^2 - 2 b K_0\pi < 0$ that corresponds to a stable focus:
\begin{align}
&\omega_l = \displaystyle -\frac{a K_0 \hskip0.1cm e^{\displaystyle t_0 \operatorname{Re} \lambda^s_1}}{2\pi} \left(c_1 \cos\left(t_0 \operatorname{Im} \lambda^s_1\right) + c_2 \sin\left(t_0 \operatorname{Im} \lambda^s_1\right)\right) + \nonumber \\
&+ \displaystyle \frac{e^{\displaystyle t_0 \operatorname{Re} \lambda^s_1}\sqrt{2 b K_0\pi - (a K_0)^2}}{2\pi} \left(c_2 \cos\left(t_0 \operatorname{Im} \lambda^s_1\right) - c_1 \sin\left(t_0 \operatorname{Im} \lambda^s_1\right)\right), \label{formula:LockInFocus}\\
&\text{where } t_0 = \frac{\operatorname{arctg}\left(\displaystyle -\frac{c_1}{c_2} \right)}{\operatorname{Im} \lambda^s_1}, c_1 = \displaystyle \frac{\pi}{2}, \hskip0.2cm c_2 = \displaystyle \frac{\pi\sqrt{(aK_0)^2 + 4{b}K_0({\pi}-\frac{1}{k})}}{\displaystyle 2\sqrt{2 b K_0\pi - (a K_0)^2}}, \nonumber \\
&\lambda^s_1=\displaystyle \frac{-aK_0 + i\sqrt{2 b K_0\pi - (a K_0)^2}}{\pi} \nonumber.
\end{align}
Rigorous derivation of (\ref{formula:LockInNode}), (\ref{formula:LockInDegNode}), and (\ref{formula:LockInFocus}) is given in \ref{sec:AppLockIn}.
The analytical and numerical results are compared in Fig.~\ref{ris:DiagramTriangleComparison}.
\section{Conclusion}
In the present work the model of PLL with impulse signals and active PI filter in the signal's phase space is described.
For the considered PLL the lock-in range is computed analytically
and obtained result are compared with numerical simulations.
\appendix

\section{The lock-in computation}
\label{sec:AppLockIn}
In this section equations (\ref{formula:LockInNode}), (\ref{formula:LockInDegNode}), and (\ref{formula:LockInFocus}) are rigorously derived. Consider the following relations
\begin{equation}
\begin{cases}
\dot{\theta}_\Delta = y,\\
\dot{y} = - a K_0 \dot{\varphi}(\theta_\Delta) y - b K_0 \varphi(\theta_\Delta).
\label{sys:PLLSysSawtooth}
\end{cases}
\end{equation}
Also we consider a normalized $2\pi$-periodic zigzag function
\begin{equation}
\varphi(\theta_\Delta)=
\begin{cases}
k \theta_\Delta, &\text{if $-\frac{1}{k} \leq \theta_\Delta \leq \frac{1}{k}$;}\\
-\frac{k}{{\pi}k-1}\theta_\Delta + \frac{{\pi}k}{{\pi}k-1}, &\text{if $\frac{1}{k} \leq \theta_\Delta \leq 2\pi-\frac{1}{k}$}
\end{cases}
\label{eq:CharSawtooth}
\end{equation}
for finite $k > \frac{1}{\pi}$ in the interval $\theta_\Delta \in \left[-\frac{1}{k}, 2\pi - \frac{1}{k}\right)$. For $k = \frac{2}{\pi}$ the function $\varphi(\theta_\Delta)$ is triangular and corresponds to (\ref{eq:PDCharTriangle}).

From $2\pi$-periodicity of (\ref{sys:PLLSysSawtooth}) it follows that for each interval the behavior of phase trajectories on the system phase plane is the same $$\theta_\Delta\in\left(-\frac{1}{k} + 2\pi j, -\frac{1}{k} + 2\pi (j+1) \right], \hskip0.2cm j \in \mathbb{Z}.$$
Thus, we can consider a single interval $\left(-\frac{1}{k}, 2\pi - \frac{1}{k}\right]$ of the phase plane of (\ref{sys:PLLSysSawtooth}).
% In Subsection \ref{sec:PhasePlaneAnalysis} for triangle characteristic $\varphi(\theta_\Delta)$ (for $k = \frac{2}{\pi}$) it is shown that equilibria $$\left(x^{eq}_{2j}, y^{eq}_{2j}\right) = \left(2\pi j, 0\right)$$ are stable equilibria and equilibria $$\left(x^{eq}_{2k+1}, y^{eq}_{2j+1}\right) = \left(2\pi j + \pi, 0\right)$$ are unstable saddle points for $\forall j \in \mathbb{Z}$. Same can be shown for sawtooth PD characteristic (\ref{eq:CharSawtooth}).

% For every other band $x\in\left[-\frac{1}{k} + 2\pi j, -\frac{1}{k} + 2\pi (j+1) \right)$ corresponding upper separatrix has the same value above the stable equilibrium $\left(x^{eq}_{2j}, y^{eq}_{2j}\right)$.

In the intervals inside $\left(-\frac{1}{k}, 2\pi - \frac{1}{k}\right]$, (\ref{sys:PLLSysSawtooth}) takes the form: \\
I. $-\frac{1}{k} < \theta_\Delta < \frac{1}{k}$ \\
\begin{align}
&\begin{cases}
\dot{\theta_\Delta} = y, \\
\dot{y} = -a K_0  k y - b K_0  k \theta_\Delta;
\end{cases}
\label{eq:PLLSys1}
\end{align}
II. $\frac{1}{k} < \theta_\Delta < 2\pi - \frac{1}{k}$
\begin{align}
&\begin{cases}
\dot{\theta_\Delta} = y, \\
\dot{y} = aK_0  \frac{k}{\pi k - 1}y + b K_0 \left(\frac{k}{\pi k -1} \theta_\Delta - \frac{\pi k}{\pi k -1}\right).
\end{cases}
\label{eq:PLLSys2}
\end{align}
In each interval there exists only one equilibrium:\\
I. $-\frac{1}{k} < \theta_\Delta < \frac{1}{k}$ \\
\\
\begin{tabular}{cc}
\hskip1cm$\begin{cases}
y_{\rm eq} = 0, \\
-a K_0  k y - b K_0  k \theta_{\rm eq} = 0;
\end{cases}$ &
\hskip3cm$\begin{cases}
y_{\rm eq} = 0, \\
\theta_{\rm eq} = 0;
\end{cases}$
\end{tabular} \\
\\
II. $\frac{1}{k} < \theta_\Delta < 2\pi - \frac{1}{k}$ \\
\\
\begin{tabular}{cc}
\hskip1cm$\begin{cases}
y_{\rm eq} = 0, \\
\frac{a K_0  k}{\pi k - 1}y_{\rm eq} + \frac{b K_0  k}{\pi k -1}\left( \theta_{\rm eq} - \pi\right) = 0.
\end{cases}$ &
\hskip2cm$\begin{cases}
y_{\rm eq} = 0, \\
\theta_{\rm eq} = \pi.
\end{cases}$
\end{tabular}

To define a type of the equilibria points, we compute the corresponding characteristic polynomial and eigenvalues. For the first equilibrium $(\theta_{\rm eq},y_{\rm eq}) = \left(0, 0\right)$ the characteristic polynomial is as follows
\begin{equation}
\chi(\lambda) = \left| \begin{array}{cc}
-\lambda & 1 \\
-b K_0  k & -a K_0  k - \lambda \end{array} \right| = \lambda^2 + aK_0 k\lambda + {b}K_0 k. \nonumber
\end{equation}
The eigenvalues of the equilibrium $(\theta_{\rm eq},y_{\rm eq}) = \left(0, 0\right)$ depend on a sign of $(a K_0 )^2 - \frac{4 b K_0 }{k}$. Here, there exist three cases: \\
\textbf{A.} $(a K_0 )^2 - \frac{4 b K_0 }{k} > 0$:
\begin{equation}
\lambda^s_{1,2} = \displaystyle \frac{-aK_0 k \pm \sqrt{(aK_0 k)^2 - 4{b}K_0 k}}{2}, \nonumber
\end{equation}
the equilibrium $\left(0, 0\right)$ is a stable node. \\
\textbf{B.} $(a K_0 )^2 - \frac{4 b K_0 }{k} = 0$:
\begin{equation}
\lambda^s_{1} = \lambda^s_{2} = \displaystyle \frac{-aK_0 k}{2}, \nonumber
\end{equation}
the equilibrium $\left(0, 0\right)$ is a stable degenerated node, or stable proper node. \\
\textbf{C.} $(a K_0 )^2 - \frac{4 b K_0 }{k} < 0$:
\begin{equation}
\lambda^s_{1,2} = \displaystyle \frac{-aK_0 k \pm i\sqrt{4{b}K_0 k - (aK_0 k)^2}}{2}, \nonumber
\end{equation}
the equilibrium $\left(0, 0\right)$ is a stable focus.  \\
Denote $(\theta^s_{\rm eq},y_{\rm eq}) = (0, 0)$.

For the second equilibrium $(\theta_{\rm eq},y_{\rm eq}) = \left(\pi, 0\right)$ we have
\begingroup
\addtolength{\jot}{1em}
\begin{align*}
&\chi(\lambda) = \left| \begin{array}{cc}
-\lambda & 1 \\
\frac{b K_0  k}{\pi k-1} & \frac{a K_0  k}{\pi k-1} - \lambda \end{array} \right| = \lambda^2 - \frac{a K_0  k}{\pi k-1} \lambda - \frac{b K_0  k}{\pi k-1}; \\
&\lambda^u_{1,2} =\displaystyle \frac{\frac{aK_0 k}{{\pi}k-1} \pm \sqrt{\left(\frac{aK_0 k}{{\pi}k-1}\right)^2 + \frac{4{b}K_0 k}{{\pi}k-1}}}{2},
\end{align*}
\endgroup
which means that $\left(\pi, 0\right)$ is always an unstable saddle for the considered parameters of the PLL. Denote $(\theta^u_{\rm eq},y_{\rm eq}) = (\pi, 0)$.

The calculation of $S^\prime(\theta^s_{\rm eq})$ from formula (\ref{eq:LockInEquivSys}) for lock-in range is in some stages. First, find two-dimensional eigenvectors $X^u_{1}$, $X^u_{2}$ of saddle point $\left(\theta^u_{\rm eq},y_{\rm eq}\right)$ from the interval $\theta_\Delta \in \left(\frac{1}{k}, 2\pi - \frac{1}{k}\right)$. Next, compute $S^\prime(\frac{1}{k})$, which is possible due to the continuity of (\ref{sys:PLLSysSawtooth}). Find two-dimensional eigenvectors $X^s_{1}$, $X^s_{2}$ of stable equilibrium $\left(\theta^s_{\rm eq},y_{\rm eq}\right)$ in the interval $\theta_\Delta \in \left(-\frac{1}{k}, \frac{1}{k}\right)$. Find a general solution of (\ref{sys:PLLSysSawtooth}) in the interval $\theta_\Delta \in \left(-\frac{1}{k}, \frac{1}{k}\right)$. Using the obtained $S^\prime(\frac{1}{k})$ as the initial data of the Cauchy problem, we can compute $S^\prime(\theta^s_{\rm eq})$.

Let us find the eigenvectors $X^u_{1}$, $X^u_{2}$ of a saddle point $\left(\theta^u_{\rm eq},y_{\rm eq}\right)$. First, find the eigenvector $X^u_{1}$:
\begingroup
\addtolength{\jot}{1em}
\begin{align*}
&\left( \begin{array}{cc}
-\lambda^u_1 & 1 \\
\frac{b K_0  k}{\pi k-1} & \frac{a K_0  k}{\pi k-1} - \lambda^u_1 \end{array} \right) X^u_1 =\mathbb{O},
\end{align*}
\begin{align*}
&\left( \begin{array}{cc}
-\displaystyle \frac{\frac{aK_0 k}{{\pi}k-1} + \sqrt{\left(\frac{aK_0 k}{{\pi}k-1}\right)^2 + \frac{4{b}K_0 k}{{\pi}k-1}}}{2} & 1 \\
\frac{b K_0  k}{\pi k-1} & \frac{a K_0  k}{\pi k-1} - \displaystyle \frac{\frac{aK_0 k}{{\pi}k-1} + \sqrt{\left(\frac{aK_0 k}{{\pi}k-1}\right)^2 + \frac{4{b}K_0 k}{{\pi}k-1}}}{2} \end{array} \right) X^u_1 =\mathbb{O},
\end{align*}
\begin{align}
&\left( \begin{array}{cc}
-\displaystyle \frac{\frac{aK_0 k}{{\pi}k-1} + \sqrt{\left(\frac{aK_0 k}{{\pi}k-1}\right)^2 + \frac{4{b}K_0 k}{{\pi}k-1}}}{2} & 1 \\
\frac{b K_0  k}{\pi k-1} & - \displaystyle \frac{\frac{aK_0 k}{{\pi}k-1} - \sqrt{\left(\frac{aK_0 k}{{\pi}k-1}\right)^2 + \frac{4{b}K_0 k}{{\pi}k-1}}}{2} \end{array} \right) X^u_1 =\mathbb{O}. \label{der:eigVecU1}
\end{align}
\endgroup
Multiply the second row of (\ref{der:eigVecU1}) by $\displaystyle \frac{\frac{aK_0 k}{{\pi}k-1} + \sqrt{\left(\frac{aK_0 k}{{\pi}k-1}\right)^2 + \frac{4{b}K_0 k}{{\pi}k-1}}}{2}$ and divide it by $\displaystyle \frac{b K_0  k}{\pi k-1}$. Then we have
\begingroup
\addtolength{\jot}{1em}
\begin{align*}
&\left( \begin{array}{cc}
-\displaystyle \frac{\frac{aK_0 k}{{\pi}k-1} + \sqrt{\left(\frac{aK_0 k}{{\pi}k-1}\right)^2 + \frac{4{b}K_0 k}{{\pi}k-1}}}{2} & 1 \\
\displaystyle \frac{\frac{aK_0 k}{{\pi}k-1} + \sqrt{\left(\frac{aK_0 k}{{\pi}k-1}\right)^2 + \frac{4{b}K_0 k}{{\pi}k-1}}}{2} & - \displaystyle \frac{\left(\left(\frac{aK_0 k}{{\pi}k-1}\right)^2 - (\frac{aK_0 k}{{\pi}k-1})^2 + \frac{4{b}K_0 k}{{\pi}k-1}\right)\left({\pi}k-1\right)}{4 b K_0  k} \end{array} \right) X^u_1 =\mathbb{O},
\end{align*}
\begin{align*}
&\left( \begin{array}{cc}
-\displaystyle \frac{\frac{aK_0 k}{{\pi}k-1} + \sqrt{\left(\frac{aK_0 k}{{\pi}k-1}\right)^2 + \frac{4{b}K_0 k}{{\pi}k-1}}}{2} & 1 \\
\displaystyle \frac{\frac{aK_0 k}{{\pi}k-1} + \sqrt{\left(\frac{aK_0 k}{{\pi}k-1}\right)^2 + \frac{4{b}K_0 k}{{\pi}k-1}}}{2} & - 1 \end{array} \right) X^u_1 =\mathbb{O}.
\end{align*}
\endgroup
Hence,
\begin{equation}
X^u_1 = \left(\begin{array}{c}
c\\
\displaystyle c\frac{\sqrt{(aK_0 k)^2 + 4{b}K_0 k({\pi}k-1)} + aK_0 k}{2({\pi}k-1)} \end{array}\right). \nonumber
\end{equation}
Let us choose $c = \displaystyle \frac{\sqrt{(aK_0 k)^2 + 4{b}K_0 k({\pi}k-1)} - aK_0 k}{2 b K_0  k}$. Then \begingroup
\addtolength{\jot}{1em}
\begin{align*}
&X^u_1 = \left(\begin{array}{c}
\displaystyle \frac{\sqrt{(aK_0 k)^2 + 4{b}K_0 k({\pi}k-1)} - aK_0 k}{2 b K_0  k} \\
\displaystyle \frac{(aK_0 k)^2 + 4{b}K_0 k({\pi}k-1) - (aK_0 k)^2}{4 b K_0  k ({\pi}k-1)}
 \end{array}\right),
\end{align*}
\begin{align*}
&X^u_1 = \left(\begin{array}{c}
\displaystyle \frac{\sqrt{(aK_0 k)^2 + 4{b}K_0 k({\pi}k-1)} - aK_0 k}{2 b K_0  k} \\
1
\end{array}\right).
\end{align*}
\endgroup

Next, find the second eigenvector $X^u_2$ in the same way:
\begingroup
\addtolength{\jot}{1em}
\begin{align*}
&\left( \begin{array}{cc}
-\lambda^u_2 & 1 \\
\frac{b K_0  k}{\pi k-1} & \frac{a K_0  k}{\pi k-1} - \lambda^u_2 \end{array} \right) X^u_2 =\mathbb{O},
\end{align*}
\begin{align*}
&\left( \begin{array}{cc}
\displaystyle \frac{\sqrt{\left(\frac{aK_0 k}{{\pi}k-1}\right)^2 + \frac{4{b}K_0 k}{{\pi}k-1}} - \frac{aK_0 k}{{\pi}k-1}}{2} & 1 \\
\frac{b K_0  k}{\pi k-1} & \frac{a K_0  k}{\pi k-1} + \displaystyle \frac{\sqrt{\left(\frac{aK_0 k}{{\pi}k-1}\right)^2 + \frac{4{b}K_0 k}{{\pi}k-1}} - \frac{aK_0 k}{{\pi}k-1}}{2} \end{array} \right) X^u_2 =\mathbb{O},
\end{align*}
\begin{align}
&\left( \begin{array}{cc}
\displaystyle \frac{\sqrt{\left(\frac{aK_0 k}{{\pi}k-1}\right)^2 + \frac{4{b}K_0 k}{{\pi}k-1}} - \frac{aK_0 k}{{\pi}k-1}}{2} & 1 \\
\frac{b K_0  k}{\pi k-1} & \displaystyle \frac{\sqrt{\left(\frac{aK_0 k}{{\pi}k-1}\right)^2 + \frac{4{b}K_0 k}{{\pi}k-1}} + \frac{aK_0 k}{{\pi}k-1}}{2} \end{array} \right) X^u_2 =\mathbb{O}. \label{der:eigVecU2}
\end{align}
\endgroup
Multiply the second row of (\ref{der:eigVecU2}) by $\displaystyle \frac{\sqrt{\left(\frac{aK_0 k}{{\pi}k-1}\right)^2 + \frac{4{b}K_0 k}{{\pi}k-1}} - \frac{aK_0 k}{{\pi}k-1}}{2}$, and divide it by $\displaystyle \frac{b K_0  k}{\pi k-1}$. Then
\begingroup
\addtolength{\jot}{1em}
\begin{align*}
&\left( \begin{array}{cc}
\displaystyle \frac{\sqrt{\left(\frac{aK_0 k}{{\pi}k-1}\right)^2 + \frac{4{b}K_0 k}{{\pi}k-1}} - \frac{aK_0 k}{{\pi}k-1}}{2} & 1 \\
\displaystyle \frac{\sqrt{\left(\frac{aK_0 k}{{\pi}k-1}\right)^2 + \frac{4{b}K_0 k}{{\pi}k-1}} - \frac{aK_0 k}{{\pi}k-1}}{2} & \displaystyle \frac{\left(\left(\frac{aK_0 k}{{\pi}k-1}\right)^2 + \frac{4{b}K_0 k}{{\pi}k-1} - \left(\frac{aK_0 k}{{\pi}k-1}\right)^2\right)\left({\pi}k-1\right)}{4 b K_0  k} \end{array} \right) X^u_2 =\mathbb{O},
\end{align*}
\begin{align*}
&\left( \begin{array}{cc}
\displaystyle \frac{\sqrt{\left(\frac{aK_0 k}{{\pi}k-1}\right)^2 + \frac{4{b}K_0 k}{{\pi}k-1}} - \frac{aK_0 k}{{\pi}k-1}}{2} & 1 \\
\displaystyle \frac{\sqrt{\left(\frac{aK_0 k}{{\pi}k-1}\right)^2 + \frac{4{b}K_0 k}{{\pi}k-1}} - \frac{aK_0 k}{{\pi}k-1}}{2} & 1 \end{array} \right) X^u_2 =\mathbb{O}.
\end{align*}
\endgroup
Hence,
\begin{equation}
X^u_2 = \left(\begin{array}{c}
-c \\
\displaystyle c\frac{\sqrt{\left(\frac{aK_0 k}{{\pi}k-1}\right)^2 + \frac{4{b}K_0 k}{{\pi}k-1}} - \frac{aK_0 k}{{\pi}k-1}}{2} \end{array}\right). \nonumber
\end{equation}
Choose $c = \displaystyle \frac{\sqrt{(aK_0 k)^2 + 4{b}K_0 k({\pi}k-1)} + aK_0 k}{2 b K_0  k}$:
\begingroup
\addtolength{\jot}{1em}
\begin{align*}
&X^u_2 = \left(\begin{array}{c}
\displaystyle -\frac{\sqrt{(aK_0 k)^2 + 4{b}K_0 k({\pi}k-1)} + aK_0 k}{2 b K_0  k} \\
\displaystyle \frac{(aK_0 k)^2 + 4{b}K_0 k({\pi}k-1) - (aK_0 k)^2}{4 b K_0  k ({\pi}k-1)}
 \end{array}\right), \nonumber \\
&X^u_2 = \left(\begin{array}{c}
\displaystyle -\frac{\sqrt{(aK_0 k)^2 + 4{b}K_0 k({\pi}k-1)} + aK_0 k}{2 b K_0  k} \\
1
 \end{array}\right).
\end{align*}
\endgroup

We can show that the direction of separatrix $S^\prime(\theta_\Delta)$ coincides with the direction of eigenvector $X^u_2$, which corresponds to eigenvalue $\lambda^u_2$. That allows us to find $S^\prime(\frac{1}{k})$. For this purpose, we write an equation of straight line, which passes through two points
\begin{align*}
&\left(x_1, y_1\right) = \left(\pi, 0\right), \nonumber \\
&\left(x_2, y_2\right) = \left(\pi - \displaystyle \frac{\sqrt{(aK_0 k)^2 + 4{b}K_0 k({\pi}k-1)} + aK_0 k}{2 b K_0  k}, 1\right).
\end{align*}
The equation takes the form
\begin{align*}
&\frac{y - 0}{1 - 0} = \frac{x - \pi}{\left(\pi - \displaystyle \frac{\sqrt{(aK_0 k)^2 + 4{b}K_0 k({\pi}k-1)} + aK_0 k}{2 b K_0  k}\right) - \pi}, \nonumber \\
&y = \frac{2 b K_0  k}{\sqrt{(aK_0 k)^2 + 4{b}K_0 k({\pi}k-1)} + aK_0 k} \left(\pi - x\right), \nonumber \\
&y = \frac{2 b K_0  k \left(\sqrt{(aK_0 k)^2 + 4{b}K_0 k({\pi}k-1)} - aK_0 k\right)}{(aK_0 k)^2 + 4{b}K_0 k({\pi}k-1) - (aK_0 k)^2} \left(\pi - x\right), \nonumber \\
&y = \frac{\sqrt{(aK_0 k)^2 + 4{b}K_0 k({\pi}k-1)} - aK_0 k}{2({\pi}k-1)} \left(\pi - x\right).
\end{align*}
Then
\begin{align*}
&S^\prime(\frac{1}{k}) = \frac{\sqrt{(aK_0 k)^2 + 4{b}K_0 k({\pi}k-1)} - aK_0 k}{2({\pi}k-1)} \left(\pi - \frac{1}{k}\right) = \\
&= \frac{\sqrt{(aK_0 )^2 + 4{b}K_0 ({\pi}-\frac{1}{k})} - aK_0 }{2}.
\end{align*}

Next, we need to find the eigenvectors of equilibrium $(\theta^s_{\rm eq},y_{\rm eq})$ and a general solution of (\ref{sys:PLLSysSawtooth}) in the interval $\left(-\frac{1}{k}, \frac{1}{k}\right)$.
It was shown that for a stable equilibrium $(\theta^s_{\rm eq},y_{\rm eq})$ in the interval $\left(-\frac{1}{k}, \frac{1}{k}\right)$ there exist three different cases, which depend on a sign of $(a K_0 )^2 - \frac{4b K_0 }{k}$. The eigenvectors $X^s_1$ and $X^s_2$ are computed in the case of stable focus only. For other cases the computation of $X^s_1$, $X^s_2$ is similar to that, considered in \ref{subsec:StableNode}.
\subsection{Stable node}
\label{subsec:StableNode}
This case corresponds to $(a K_0 )^2 - \frac{4 b K_0 }{k} > 0$. Let us find the eigenvectors $X^s_1$, $X^s_2$:
\begingroup
\addtolength{\jot}{1em}
\begin{align*}
&\left( \begin{array}{cc}
-\lambda^s_1 & 1 \\
-b K_0  k & -a K_0  k - \lambda^s_1 \end{array} \right) X^s_1 =\mathbb{O},
\end{align*}
\begin{align*}
&\left( \begin{array}{cc}
\displaystyle \frac{aK_0 k - \sqrt{(aK_0 k)^2 - 4{b}K_0 k}}{2} & 1 \\
-b K_0  k & \displaystyle -aK_0 k +\frac{aK_0 k - \sqrt{(aK_0 k)^2 - 4{b}K_0 k}}{2} \end{array} \right) X^s_1 =\mathbb{O},
\end{align*}
\begin{align}
&\left( \begin{array}{cc}
\displaystyle \frac{aK_0 k - \sqrt{(aK_0 k)^2 - 4{b}K_0 k}}{2} & 1 \\
-b K_0  k & \displaystyle -\frac{aK_0 k + \sqrt{(aK_0 k)^2 - 4{b}K_0 k}}{2} \end{array} \right) X^s_1 =\mathbb{O}. \label{der:eigVecS1}
\end{align}
\endgroup
Multiply the second row of (\ref{der:eigVecS1}) by $\displaystyle \frac{aK_0 k - \sqrt{(aK_0 k)^2 - 4{b}K_0 k}}{2}$, and divide it by $b K_0  k$:
\begingroup
\addtolength{\jot}{1em}
\begin{align*}
&\left( \begin{array}{cc}
\displaystyle \frac{aK_0 k - \sqrt{(aK_0 k)^2 - 4{b}K_0 k}}{2} & 1 \\
-\displaystyle \frac{aK_0 k - \sqrt{(aK_0 k)^2 - 4{b}K_0 k}}{2} & \displaystyle -\frac{\left(aK_0 k\right)^2 - (aK_0 k)^2 + 4{b}K_0 k}{4 b K_0  k} \end{array} \right) X^s_1 =\mathbb{O},
\end{align*}
\begin{align*}
&\left( \begin{array}{cc}
\displaystyle \frac{aK_0 k - \sqrt{(aK_0 k)^2 - 4{b}K_0 k}}{2} & 1 \\
-\displaystyle \frac{aK_0 k - \sqrt{(aK_0 k)^2 - 4{b}K_0 k}}{2} & -1 \end{array} \right) X^s_1 =\mathbb{O}, \end{align*}
\endgroup
\begin{equation}
X^s_1 = \left(\begin{array}{c}
- c \\
c \displaystyle \frac{aK_0 k - \sqrt{(aK_0 k)^2 - 4{b}K_0 k}}{2} \end{array}\right). \nonumber
\end{equation}
Choose $c = -1$. Then
\begin{equation}
X^s_1 = \left(\begin{array}{c}
1 \\
\displaystyle \frac{\sqrt{(aK_0 k)^2 - 4{b}K_0 k} - aK_0 k}{2} \end{array}\right). \nonumber
\end{equation}
Next, find eigenvector $X^s_2$:
\begingroup
\addtolength{\jot}{1em}
\begin{align*}
&\left( \begin{array}{cc}
-\lambda^s_2 & 1 \\
-b K_0  k & -a K_0  k - \lambda^s_2 \end{array} \right) X^s_2 =\mathbb{O},
\end{align*}
\begin{align*}
&\left( \begin{array}{cc}
\displaystyle \frac{aK_0 k + \sqrt{(aK_0 k)^2 - 4{b}K_0 k}}{2} & 1 \\
-b K_0  k & \displaystyle -aK_0 k +\frac{aK_0 k + \sqrt{(aK_0 k)^2 - 4{b}K_0 k}}{2} \end{array} \right) X^s_2 =\mathbb{O},
\end{align*}
\begin{align}
&\left( \begin{array}{cc}
\displaystyle \frac{aK_0 k + \sqrt{(aK_0 k)^2 - 4{b}K_0 k}}{2} & 1 \\
-b K_0  k & \displaystyle \frac{\sqrt{(aK_0 k)^2 - 4{b}K_0 k} - aK_0 k}{2} \end{array} \right) X^s_2 =\mathbb{O}. \label{der:eigVecS2}
\end{align}
\endgroup
Multiply the second row of (\ref{der:eigVecS2}) by $\displaystyle \frac{aK_0 k + \sqrt{(aK_0 k)^2 - 4{b}K_0 k}}{2}$, and divide it by $b K_0  k$:
\begingroup
\addtolength{\jot}{1em}
\begin{align*}
&\left( \begin{array}{cc}
\displaystyle \frac{aK_0 k + \sqrt{(aK_0 k)^2 - 4{b}K_0 k}}{2} & 1 \\
-\displaystyle \frac{aK_0 k + \sqrt{(aK_0 k)^2 - 4{b}K_0 k}}{2} & \displaystyle \frac{\left(aK_0 k\right)^2 - 4{b}K_0 k - (aK_0 k)^2}{4 b K_0  k} \end{array} \right) X^s_2 =\mathbb{O},
\end{align*}
\begin{align*}
&\left( \begin{array}{cc}
\displaystyle \frac{aK_0 k + \sqrt{(aK_0 k)^2 - 4{b}K_0 k}}{2} & 1 \\
-\displaystyle \frac{aK_0 k + \sqrt{(aK_0 k)^2 - 4{b}K_0 k}}{2} & -1 \end{array} \right) X^s_2 =\mathbb{O}, \end{align*}
\endgroup
\begin{equation}
X^s_2 = \left(\begin{array}{c}
- c \\
c \displaystyle \frac{aK_0 k + \sqrt{(aK_0 k)^2 - 4{b}K_0 k}}{2} \end{array}\right). \nonumber
\end{equation}
Choose $c=-1$. Then
\begin{equation}
X^s_2 = \left(\begin{array}{c}
1 \\
\displaystyle -\frac{aK_0 k + \sqrt{(aK_0 k)^2 - 4{b}K_0 k}}{2} \end{array}\right). \nonumber
\end{equation}
In the interval $\theta_\Delta \in \left(-\frac{1}{k}, \frac{1}{k}\right)$ for $\left(\theta^s_{\rm eq},y_{\rm eq}\right) = \left(0, 0\right)$ being a node, a general solution of (\ref{sys:PLLSysSawtooth}) has the form:
\begin{align}
\begin{cases}
\theta_\Delta(t) = c_1{\hskip0.1cm}e^{\displaystyle \lambda^s_1 t} + c_2{\hskip0.1cm}e^{\displaystyle \lambda^s_2 t}, \\
\\
y(t) = \displaystyle -c_1\frac{aK_0 k - \sqrt{(aK_0 k)^2 - 4{b}K_0 k}}{2}{\hskip0.1cm}e^{\displaystyle \lambda^s_1 t} - \displaystyle c_2\frac{aK_0 k + \sqrt{(aK_0 k)^2 - 4{b}K_0 k}}{2}{\hskip0.1cm}e^{\displaystyle \lambda^s_2 t}.
\end{cases} \label{eq:CommSolNode}
\end{align}
Let us find coefficients $c_1$, $c_2$ of (\ref{eq:CommSolNode}) for the solution of the Cauchy problem with initial conditions $\theta_\Delta(0) = \frac{1}{k}$, $y(0) = \frac{\sqrt{(aK_0 )^2 + 4{b}K_0 ({\pi}-\frac{1}{k})} - aK_0 }{2}$, which coincide with $S^\prime(\frac{1}{k})$.\\
At moment $t=0$ we have
\begin{equation}
\begin{cases}
\displaystyle \frac{1}{k} = c_1 + c_2, \\
\\
\displaystyle \frac{\sqrt{(aK_0 )^2 + 4{b}K_0 ({\pi}-\frac{1}{k})} - aK_0 }{2} = \\
\\
\hskip1.5cm  = \displaystyle - c_1\frac{aK_0 k - \sqrt{(aK_0 k)^2 - 4{b}K_0 k}}{2} - \displaystyle c_2\frac{aK_0 k + \sqrt{(aK_0 k)^2 - 4{b}K_0 k}}{2},
\end{cases} \nonumber
\end{equation}
\begin{equation}
\begin{cases}
\displaystyle c_2 = \frac{1}{k} - c_1, \\
\\
\displaystyle \frac{\sqrt{(aK_0 )^2 + 4{b}K_0 ({\pi}-\frac{1}{k})} - aK_0 }{2} + \frac{aK_0 k + \sqrt{(aK_0 k)^2 - 4{b}K_0 k}}{2k}= \\
=\displaystyle -c_1\frac{aK_0 k - \sqrt{(aK_0 k)^2 - 4{b}K_0 k}}{2} + \displaystyle c_1\frac{aK_0 k + \sqrt{(aK_0 k)^2 - 4{b}K_0 k}}{2},
\end{cases} \nonumber
\end{equation}
\begin{equation}
\begin{cases}
\displaystyle c_2 = \frac{1}{k} - c_1, \\
\\
\displaystyle \frac{\sqrt{(aK_0 )^2 + 4{b}K_0 ({\pi}-\frac{1}{k})}}{2} + \frac{\sqrt{(aK_0 )^2 - \frac{4{b}K_0 }{k}}}{2} = \displaystyle c_1 k \sqrt{(aK_0 )^2 - \frac{4{b}K_0 }{k}},
\end{cases}  \nonumber
\end{equation}
\begin{equation}
\begin{cases}
\displaystyle c_2 = \frac{1}{k} - c_1, \\
\\
\displaystyle c_1  = \left(\displaystyle \frac{\sqrt{(aK_0 )^2 + 4{b}K_0 ({\pi}-\frac{1}{k})}}{\sqrt{(aK_0 )^2 - \frac{4{b}K_0 }{k}}} + 1\right):2k,
\end{cases} \nonumber
\end{equation}
\begin{equation}
\begin{cases}
\displaystyle c_1  = \left(\displaystyle \frac{\sqrt{(aK_0 )^2 + 4{b}K_0 ({\pi}-\frac{1}{k})}}{\sqrt{(aK_0 )^2 - \frac{4{b}K_0 }{k}}} + 1\right):2k, \\
c_2 = \left(1 - \displaystyle \frac{\sqrt{(aK_0 )^2 + 4{b}K_0 ({\pi}-\frac{1}{k})}}{\sqrt{(aK_0 )^2 - \frac{4{b}K_0 }{k}}}\right):2k.
\end{cases}
\label{rel:Const12Node}
\end{equation}
Finally, find $y(t_0)$ under the condition $\theta_\Delta(t_0) = 0$. The value of $y(t_0)$ corresponds to $S^\prime(\theta^s_{\rm eq})$. For this purpose, we express $y(t_0)$ in terms of $c_1$, $c_2$ from (\ref{rel:Const12Node}). Then
\begin{equation}
\begin{cases}
0 = c_1{\hskip0.1cm}e^{\displaystyle \lambda^s_1 t_0} + c_2{\hskip0.1cm}e^{\displaystyle \lambda^s_2 t_0}, \\
\\
y(t_0) = \displaystyle -c_1\frac{aK_0 k - \sqrt{(aK_0 k)^2 - 4{b}K_0 k}}{2}{\hskip0.1cm}e^{\displaystyle \lambda^s_1 t_0} - \displaystyle c_2\frac{aK_0 k + \sqrt{(aK_0 k)^2 - 4{b}K_0 k}}{2}{\hskip0.1cm}e^{\displaystyle \lambda^s_2 t_0},
\end{cases} \nonumber
\end{equation}

\begin{equation}
\begin{cases}
\displaystyle -\frac{c_1}{c_2}  = e^{\displaystyle \left(\lambda^s_2 - \lambda^s_1\right) t_0}, \\
\\
y(t_0) = \displaystyle -c_1\frac{aK_0 k - \sqrt{(aK_0 k)^2 - 4{b}K_0 k}}{2}{\hskip0.1cm}e^{\displaystyle \lambda^s_1 t_0} - \displaystyle c_2\frac{aK_0 k + \sqrt{(aK_0 k)^2 - 4{b}K_0 k}}{2}{\hskip0.1cm}e^{\displaystyle \lambda^s_2 t_0},
\end{cases} \nonumber
\end{equation}

\begin{equation}
\begin{cases}
\displaystyle -\frac{c_1}{c_2}  = e^{\displaystyle \left(-\frac{\sqrt{(aK_0 k)^2 - 4{b}K_0 k} + aK_0 k}{2} - \frac{\sqrt{(aK_0 k)^2 - 4{b}K_0 k} - aK_0 k}{2}\right) t_0}, \\
\\
y(t_0) = \displaystyle -c_1\frac{aK_0 k - \sqrt{(aK_0 k)^2 - 4{b}K_0 k}}{2}{\hskip0.1cm}e^{\displaystyle \lambda^s_1 t_0} - \displaystyle c_2\frac{aK_0 k + \sqrt{(aK_0 k)^2 - 4{b}K_0 k}}{2}{\hskip0.1cm}e^{\displaystyle \lambda^s_2 t_0},
\end{cases} \nonumber
\end{equation}

\begin{equation}
\begin{cases}
\displaystyle -\frac{c_1}{c_2}  = e^{\displaystyle \left(-\sqrt{(aK_0 k)^2 - 4{b}K_0 k}\right) t_0}, \\
\\
y(t_0) = \displaystyle -c_1\frac{aK_0 k - \sqrt{(aK_0 k)^2 - 4{b}K_0 k}}{2}{\hskip0.1cm}e^{\displaystyle \lambda^s_1 t_0} - \displaystyle c_2\frac{aK_0 k + \sqrt{(aK_0 k)^2 - 4{b}K_0 k}}{2}{\hskip0.1cm}e^{\displaystyle \lambda^s_2 t_0},
\end{cases} \nonumber
\end{equation}

% \begin{equation}
% \begin{cases}
% \displaystyle -\frac{\left(\displaystyle \frac{\sqrt{(aK_0 )^2 + 4{b}K_0 ({\pi}-\frac{1}{k})}}{\sqrt{(aK_0 )^2 - \frac{4{b}K_0 }{k}}} + 1\right)}{\left(1-\displaystyle \frac{\sqrt{(aK_0 )^2 + 4{b}K_0 ({\pi}-\frac{1}{k})}}{\sqrt{(aK_0 )^2 - \frac{4{b}K_0 }{k}}}\right)}  = e^{\displaystyle \left(-\sqrt{(aK_0 k)^2 - 4{b}K_0 k}\right) t_0}, \\
% \\
% y(t_0) = \displaystyle -c_1\frac{aK_0 k - \sqrt{(aK_0 k)^2 - 4{b}K_0 k}}{2}{\hskip0.1cm}e^{\displaystyle \lambda^s_1 t_0} - \displaystyle c_2\frac{aK_0 k + \sqrt{(aK_0 k)^2 - 4{b}K_0 k}}{2}{\hskip0.1cm}e^{\displaystyle \lambda^s_2 t_0}.
% \end{cases} \nonumber
% \end{equation}

% \begin{equation}
% \begin{cases}
% \displaystyle -\frac{\left(\displaystyle \sqrt{(aK_0 )^2 + 4{b}K_0 ({\pi}-\frac{1}{k})} + \sqrt{(aK_0 )^2 - \frac{4{b}K_0 }{k}}\right)}{\left(1-\displaystyle \frac{\sqrt{(aK_0 )^2 + 4{b}K_0 ({\pi}-\frac{1}{k})}}{\sqrt{(aK_0 )^2 - \frac{4{b}K_0 }{k}}}\right)}  = e^{\displaystyle \left(-\sqrt{(aK_0 k)^2 - 4{b}K_0 k}\right) t_0}, \\
% \\
% y(t_0) = \displaystyle -c_1\frac{aK_0 k - \sqrt{(aK_0 k)^2 - 4{b}K_0 k}}{2}{\hskip0.1cm}e^{\displaystyle \lambda^s_1 t_0} - \displaystyle c_2\frac{aK_0 k + \sqrt{(aK_0 k)^2 - 4{b}K_0 k}}{2}{\hskip0.1cm}e^{\displaystyle \lambda^s_2 t_0}.
% \end{cases} \nonumber
% \end{equation}

\begin{equation}
\begin{cases}
\ln \left(\displaystyle -\frac{c_1}{c_2}\right)  = -\left(\sqrt{(aK_0 k)^2 - 4{b}K_0 k}\right) t_0, \\
\\
y(t_0) = \displaystyle -c_1\frac{aK_0 k - \sqrt{(aK_0 k)^2 - 4{b}K_0 k}}{2}{\hskip0.1cm}e^{\displaystyle \lambda^s_1 t_0} - \displaystyle c_2\frac{aK_0 k + \sqrt{(aK_0 k)^2 - 4{b}K_0 k}}{2}{\hskip0.1cm}e^{\displaystyle \lambda^s_2 t_0},
\end{cases} \nonumber
\end{equation}

\begin{equation}
\begin{cases}
t_0 = \frac{\displaystyle \ln \left( -\frac{c_2}{c_1}\right)}{ \displaystyle \sqrt{(aK_0 k)^2 - 4{b}K_0 k}}, \\
\\
y(t_0) = \displaystyle -c_1\frac{aK_0 k - \sqrt{(aK_0 k)^2 - 4{b}K_0 k}}{2}{\hskip0.1cm}e^{\displaystyle \lambda^s_1 t_0} - \displaystyle c_2\frac{aK_0 k + \sqrt{(aK_0 k)^2 - 4{b}K_0 k}}{2}{\hskip0.1cm}e^{\displaystyle \lambda^s_2 t_0},
\end{cases} \nonumber
\end{equation}
Transform the following expression
\begin{align*}
&e^{\displaystyle \lambda^s_1 t_0} = e^{\left(\displaystyle \ln \left( -\frac{c_2}{c_1}\right) \frac{\lambda^s_1}{ \displaystyle \sqrt{(aK_0 k)^2 - 4{b}K_0 k}}\right)} = e^{\left(\displaystyle \ln \left( -\frac{c_2}{c_1}\right) \frac{\sqrt{(aK_0 k)^2 - 4{b}K_0 k} - a K_0  k}{\displaystyle 2 \sqrt{(aK_0 k)^2 - 4{b}K_0 k}}\right)} = \\
& = \left(e^{\displaystyle \ln \left( -\frac{c_2}{c_1}\right)} \right)^{\left(\displaystyle \frac{1}{2} - \frac{a K_0  k}{\displaystyle 2 \sqrt{(aK_0 k)^2 - 4{b}K_0 k}}\right)} = \left( -\frac{c_2}{c_1}\right)^{\left(\displaystyle \frac{1}{2} - \frac{a K_0  k}{\displaystyle 2 \sqrt{(aK_0 k)^2 - 4{b}K_0 k}}\right)}.
\end{align*}
Similarly,
\begin{align*}
&e^{\displaystyle \lambda^s_2 t_0} = e^{\left( \displaystyle \ln \left( -\frac{c_2}{c_1}\right) \frac{\lambda^s_2}{\displaystyle \sqrt{(aK_0 k)^2 - 4{b}K_0 k}}\right)} = e^{-\left(\displaystyle \ln \left( -\frac{c_2}{c_1}\right) \frac{\sqrt{(aK_0 k)^2 - 4{b}K_0 k} + a K_0  k}{\displaystyle 2 \sqrt{(aK_0 k)^2 - 4{b}K_0 k}}\right)} = \\
& = \left(e^{\displaystyle \ln \left( -\frac{c_2}{c_1}\right)} \right)^{-\left(\displaystyle \frac{1}{2} + \frac{a K_0 }{\displaystyle 2 \sqrt{(aK_0 )^2 - \frac{4{b}K_0 }{k}}}\right)} = \left(- \frac{c_2}{c_1}\right)^{\left(\displaystyle \frac{1}{2} - \frac{a K_0  k}{\displaystyle 2 \sqrt{(aK_0 k)^2 - 4{b}K_0 k}}\right) \displaystyle-1}.
\end{align*}
Then
\begin{align*}
&y(t_0) = \displaystyle -c_1\frac{aK_0 k - \sqrt{(aK_0 k)^2 - 4{b}K_0 k}}{2}{\hskip0.1cm}e^{\displaystyle \lambda^s_1 t_0} - \displaystyle c_2\frac{aK_0 k + \sqrt{(aK_0 k)^2 - 4{b}K_0 k}}{2}{\hskip0.1cm}e^{\displaystyle \lambda^s_2 t_0},
\end{align*}
\begin{align*}
&y(t_0) = \displaystyle -c_1\frac{aK_0 k - \sqrt{(aK_0 k)^2 - 4{b}K_0 k}}{2}\left( -\frac{c_2}{c_1}\right)^{\left(\displaystyle \frac{1}{2} - \frac{a K_0  k}{\displaystyle 2 \sqrt{(aK_0 k)^2 - 4{b}K_0 k}}\right)} - \\
&- \displaystyle c_2\frac{aK_0 k + \sqrt{(aK_0 k)^2 - 4{b}K_0 k}}{2}{\hskip0.1cm}\left(- \frac{c_2}{c_1}\right)^{\left(\displaystyle \frac{1}{2} - \frac{a K_0  k}{\displaystyle 2 \sqrt{(aK_0 k)^2 - 4{b}K_0 k}}\right) \displaystyle-1},
\end{align*}
\begin{align*}
&y(t_0) = \displaystyle -c_1\frac{aK_0 k - \sqrt{(aK_0 k)^2 - 4{b}K_0 k}}{2}\left( -\frac{c_2}{c_1}\right)^{\left(\displaystyle \frac{1}{2} - \frac{a K_0  k}{\displaystyle 2 \sqrt{(aK_0 k)^2 - 4{b}K_0 k}}\right)} + \\
&+ \displaystyle c_1\frac{aK_0 k + \sqrt{(aK_0 k)^2 - 4{b}K_0 k}}{2}{\hskip0.1cm}\left(- \frac{c_2}{c_1}\right)^{\left(\displaystyle \frac{1}{2} - \frac{a K_0  k}{\displaystyle 2 \sqrt{(aK_0 k)^2 - 4{b}K_0 k}}\right)}.
\end{align*}

As a result, for the case $(aK_0 k)^2 - 4{b}K_0 k > 0$, when a stable equilibrium $\left(\theta^s_{\rm eq},y_{\rm eq}\right)$ is a stable node, $S^\prime(\theta^s_{\rm eq})$ can be found from the following formula
\begin{align}
&S^\prime(\theta^s_{\rm eq}) = \displaystyle c_1 \sqrt{(aK_0 k)^2 - 4{b}K_0 k}\left( -\frac{c_2}{c_1}\right)^{\left(\displaystyle \frac{1}{2} - \frac{a K_0  k}{\displaystyle 2 \sqrt{(aK_0 k)^2 - 4{b}K_0 k}}\right)}, \label{rel:S0Node}
\end{align}
where
\begin{equation}
\displaystyle c_1  = \left(\displaystyle \frac{\sqrt{(aK_0 )^2 + 4{b}K_0 ({\pi}-\frac{1}{k})}}{\sqrt{(aK_0 )^2 - \frac{4{b}K_0 }{k}}} + 1\right):2k, \hskip1cm c_2 = \left(1 - \displaystyle \frac{\sqrt{(aK_0 )^2 + 4{b}K_0 ({\pi}-\frac{1}{k})}}{\sqrt{(aK_0 )^2 - \frac{4{b}K_0 }{k}}}\right):2k. \nonumber
\end{equation}

\subsection{Stable focus}
\label{subsec:StableFocus}
This case corresponds to $(a K_0 )^2 - \frac{4 b K_0 }{k} < 0$. The eigenvectors $X^s_1$, $X^s_2$ are found in the same way as in \ref{subsec:StableNode}:
\begin{equation}
X^s_1 = \left(\begin{array}{c}
1 \\
\displaystyle -\frac{aK_0 k - i\sqrt{4{b}K_0 k - (aK_0 k)^2}}{2} \end{array}\right), \hskip1cm X^s_2 = \left(\begin{array}{c}
1 \\
\displaystyle -\frac{aK_0 k + i\sqrt{4{b}K_0 k - (aK_0 k)^2}}{2} \end{array}\right). \nonumber
\end{equation}
The eigenvectors $X^s_1$, $X^s_2$ can be represented as
\begin{equation}
X^s_{1,2}(t) = U^s_{1,2} + i V^s_{1,2}, \nonumber
\end{equation}
where $U^s_{1,2}$, $V^s_{1,2}$ are real two-dimensional vectors:
\begin{align*}
&U^s_1 = \left(\begin{array}{c}
1 \\
\displaystyle -\frac{aK_0 k}{2} \end{array}\right), \hskip1cm V^s_1 = \left(\begin{array}{c}
0 \\
\displaystyle \frac{\sqrt{4{b}K_0 k - (aK_0 k)^2}}{2} \end{array}\right), \\
&U^s_2 = \left(\begin{array}{c}
1 \\
\displaystyle -\frac{aK_0 k}{2} \end{array}\right), \hskip1cm V^s_2 = \left(\begin{array}{c}
0 \\
\displaystyle -\frac{\sqrt{4{b}K_0 k - (aK_0 k)^2}}{2} \end{array}\right).
\end{align*}
Consider a solution of (\ref{sys:PLLSysSawtooth}), which corresponds to the eigenvalue $\lambda^s_1$:
\begin{align*}
&Y_1(t) = e^{\lambda^s_1 t}X^s_1 = e^{\displaystyle \left(\operatorname{Re} \lambda^s_1 + i \operatorname{Im} \lambda^s_1\right) t} \hskip0.1cm \left(U^s_1 + i V^s_1\right) = \\
& = e^{\displaystyle t \operatorname{Re} \lambda^s_1} \left(\cos\left(t \operatorname{Im} \lambda^s_1 \right) + i \sin\left(t \operatorname{Im} \lambda^s_1 \right)\right) \left(U^s_1 + i V^s_1\right) = \\
&e^{\displaystyle t \operatorname{Re} \lambda^s_1} \left(U^s_1 \cos\left(t \operatorname{Im} \lambda^s_1\right) - V^s_1 \sin\left(t \operatorname{Im} \lambda^s_1 \right)\right) + \\
&+ i e^{\displaystyle t \operatorname{Re} \lambda^s_1} \left(U^s_1 \sin\left(t \operatorname{Im} \lambda^s_1 \right) + V^s_1 \cos\left(t \operatorname{Im} \lambda^s_1 \right)\right).
\end{align*}
A general solution of (\ref{sys:PLLSysSawtooth}) takes the form
\begin{align*}
&Y(t) = c_1 e^{\displaystyle t \operatorname{Re} \lambda^s_1} \left(U^s_1 \cos\left(t \operatorname{Im} \lambda^s_1 \right) - V^s_1 \sin\left(t \operatorname{Im} \lambda^s_1 \right)\right) + \\
&+ c_2 e^{\displaystyle t \operatorname{Re} \lambda^s_1} \left(U^s_1 \sin\left(t \operatorname{Im} \lambda^s_1 \right) + V^s_1 \cos\left(t \operatorname{Im} \lambda^s_1 \right)\right) = \\
&= e^{\displaystyle t \operatorname{Re} \lambda^s_1} U^s_1 \left(c_1 \cos\left(t \operatorname{Im} \lambda^s_1 \right) + c_2 \sin\left(t \operatorname{Im} \lambda^s_1 \right)\right) + \\
& + e^{\displaystyle t \operatorname{Re} \lambda^s_1} V^s_1 \left(c_2 \cos\left(t \operatorname{Im} \lambda^s_1 \right) - c_1 \sin\left(t \operatorname{Im} \lambda^s_1 \right)\right).
\end{align*}
In other words,
\begin{align}
&\begin{cases}
\theta_\Delta(t) = e^{\displaystyle t \operatorname{Re} \lambda^s_1} \left(c_1 \cos\left(t \operatorname{Im} \lambda^s_1 \right) + c_2 \sin\left(t \operatorname{Im} \lambda^s_1 \right)\right), \\
y(t) = \displaystyle -\frac{a K_0  k \hskip0.1cm e^{\displaystyle t \operatorname{Re} \lambda^s_1}}{2} \left(c_1 \cos\left(t \operatorname{Im} \lambda^s_1 \right) + c_2 \sin\left(t \operatorname{Im} \lambda^s_1 \right)\right) + \\
+ \displaystyle \frac{e^{\displaystyle t \operatorname{Re} \lambda^s_1}\sqrt{4{b}K_0 k - (aK_0 k)^2}}{2} \left(c_2 \cos\left(t \operatorname{Im} \lambda^s_1 \right) - c_1 \sin\left(t \operatorname{Im} \lambda^s_1 \right)\right).
\end{cases} \label{eq:CommSolFocus}
\end{align}
Let us find the coefficients $c_1$, $c_2$ for the solution of the Cauchy problem with the initial data $\theta_\Delta(0) = \frac{1}{k}$, $y(0) = \frac{\sqrt{(aK_0 )^2 + 4{b}K_0 ({\pi}-\frac{1}{k})} - aK_0 }{2}$, similarly to \ref{subsec:StableNode}.
\begin{align*}
&\begin{cases}
\displaystyle \frac{1}{k} = e^{\displaystyle 0 \operatorname{Re} \lambda^s_1} \left(c_1 \cos\left(0 \operatorname{Im} \lambda^s_1\right) + c_2 \sin\left(0 \operatorname{Im} \lambda^s_1\right)\right), \\
\displaystyle \frac{\sqrt{(aK_0 )^2 + 4{b}K_0 ({\pi}-\frac{1}{k})} - aK_0 }{2} = \\
= \displaystyle -\frac{a K_0  k \hskip0.1cm e^{\displaystyle 0 \operatorname{Re} \lambda^s_1}}{2} \left(c_1 \cos\left(0 \operatorname{Im} \lambda^s_1\right) + c_2 \sin\left(0 \operatorname{Im} \lambda^s_1\right)\right) + \\
+ \displaystyle \frac{e^{\displaystyle 0 \operatorname{Re} \lambda^s_1}\sqrt{4{b}K_0 k - (aK_0 k)^2}}{2} \left(c_2 \cos\left(0 \operatorname{Im} \lambda^s_1\right) - c_1 \sin\left(0 \operatorname{Im} \lambda^s_1\right)\right),
\end{cases}
\end{align*}
\begin{align*}
&\begin{cases}
\displaystyle \frac{1}{k} = c_1, \\
\displaystyle \frac{\sqrt{(aK_0 )^2 + 4{b}K_0 ({\pi}-\frac{1}{k})} - aK_0 }{2} = \displaystyle -\frac{a K_0  k }{2} c_1 + \displaystyle \frac{\sqrt{4{b}K_0 k - (aK_0 k)^2}}{2} c_2,
\end{cases}
\end{align*}
\begin{align*}
&\begin{cases}
c_1 = \displaystyle \frac{1}{k}, \\
c_2 = \displaystyle \frac{\sqrt{(aK_0 )^2 + 4{b}K_0 ({\pi}-\frac{1}{k})}}{\displaystyle k\sqrt{\frac{4{b}K_0 }{k} - (aK_0 )^2}}.
\end{cases}
\end{align*}
Next, let us find $t_0$ such that $\theta_\Delta(t_0) = 0$.
\begin{equation}
0 = e^{\displaystyle t_0 \operatorname{Re} \lambda^s_1} \left(c_1 \cos\left(t_0 \operatorname{Im} \lambda^s_1\right) + c_2 \sin\left(t_0 \operatorname{Im} \lambda^s_1\right)\right), \nonumber
\end{equation}
\begin{equation}
\displaystyle -\frac{c_1}{c_2} = \operatorname{tg} \left(t_0 \operatorname{Im} \lambda^s_1\right), \nonumber
\end{equation}
\begin{equation}
t_0 = \frac{\operatorname{arctg}\left(\displaystyle -\frac{c_1}{c_2} \right)}{\operatorname{Im} \lambda^s_1}. \nonumber
\end{equation}
Finally, all the unknowns for $y(t_0)$ from (\ref{eq:CommSolFocus}) are found and $S^\prime(\theta^s_{\rm eq})$ is as follows
\begin{align}
&S^\prime(\theta^s_{\rm eq}) = y(t_0) = \displaystyle -\frac{a K_0  k \hskip0.1cm e^{\displaystyle t_0 \operatorname{Re} \lambda^s_1}}{2} \left(c_1 \cos\left(t_0 \operatorname{Im} \lambda^s_1\right) + c_2 \sin\left(t_0 \operatorname{Im} \lambda^s_1\right)\right) + \nonumber \\
&+ \displaystyle \frac{e^{\displaystyle t_0 \operatorname{Re} \lambda^s_1}\sqrt{4{b}K_0 k - (aK_0 k)^2}}{2} \left(c_2 \cos\left(t_0 \operatorname{Im} \lambda^s_1\right) - c_1 \sin\left(t_0 \operatorname{Im} \lambda^s_1\right)\right), \label{rel:S0Focus}
\end{align}
where
\begin{align*}
&t_0 = \frac{\operatorname{arctg}\left(\displaystyle -\frac{c_1}{c_2} \right)}{\operatorname{Im} \lambda^s_1}, \hskip1cm \lambda^s_{1} = \displaystyle \frac{-aK_0 k + i\sqrt{4{b}K_0 k - (aK_0 k)^2}}{2},\\
&c_1 = \displaystyle \frac{1}{k}, \hskip1cm c_2 = \displaystyle \frac{\sqrt{(aK_0 )^2 + 4{b}K_0 ({\pi}-\frac{1}{k})}}{\displaystyle k\sqrt{\frac{4{b}K_0 }{k} - (aK_0 )^2}}. \nonumber
\end{align*}
\subsection{Stable degenerated node}
\label{subsec:StableDegNode}
This case corresponds to $(a K_0 )^2 - \frac{4 b K_0 }{k} = 0$. In this case the eigenvalues $\lambda^s_1$ and $\lambda^s_1$ coincide: $$\lambda^s := \displaystyle -\frac{a K_0  k}{2} = \lambda^s_1 = \lambda^s_2.$$
A stable equilibrium $\left(\theta^s_{\rm eq},y_{\rm eq}\right)$ is a stable degenerated node, or a stable proper node.

For the characteristic matrix
\begin{equation}
\left( \begin{array}{cc}
-\lambda^s & 1 \\
-b K_0  k & -a K_0  k - \lambda^s \end{array} \right) \nonumber
\end{equation}
of (\ref{sys:PLLSysSawtooth}) it is shown that $\left(\theta^s_{\rm eq},y_{\rm eq}\right)$ is a stable degenerated node.
Find the eigenvector $X^s$ corresponding to the eigenvalue $\lambda^s$ of algebraic multiplicity two:
\begin{equation}
\left( \begin{array}{cc}
-\lambda^s & 1 \\
-b K_0  k & -a K_0  k - \lambda^s \end{array} \right) X^s =\mathbb{O}, \nonumber
\end{equation}
\begin{equation}
\left( \begin{array}{cc}
\displaystyle \frac{aK_0 k}{2} & 1 \\
-b K_0  k & \displaystyle -\frac{aK_0 k}{2} \end{array} \right) X^s =\mathbb{O}. \nonumber
\end{equation}
Adding the first row, multiplied by $\displaystyle \frac{a K_0  k}{2}$, to the second row, we have:
\begin{equation}
\left( \begin{array}{cc}
\displaystyle \frac{aK_0 k}{2} & 1 \\
\displaystyle \frac{\left(a K_0  k\right)^2}{4}-b K_0  k & 0 \end{array} \right) X^s =\mathbb{O}. \nonumber
\end{equation}

The eigenvector $X^s$ can be written as
\begin{equation}
X^s = \left( \begin{array}{c}
c \\
\displaystyle -c\frac{a K_0  k}{2} \end{array} \right). \nonumber
\end{equation}
Choose $c = 1$. To find a general solution of (\ref{sys:PLLSysSawtooth}), we need to additionally find the first associated vector $X^s_1$:
\begin{equation}
\left( \begin{array}{cc}
\displaystyle \frac{aK_0 k}{2} & 1 \\
-b K_0  k & \displaystyle -\frac{aK_0 k}{2} \end{array} \right) X^s_1 = \left( \begin{array}{c}
1 \\
\displaystyle -\frac{a K_0  k}{2} \end{array} \right), \nonumber
\end{equation}
\begin{equation}
\left( \begin{array}{cc}
\displaystyle \frac{aK_0 k}{2} & 1 \\
\displaystyle \frac{\left(aK_0 k\right)^2}{4}-b K_0  k & 0 \end{array} \right) X^s_1 = \left( \begin{array}{c}
1 \\
0 \end{array} \right), \nonumber
\end{equation}
\begin{equation}
X^s_1 = \left( \begin{array}{c}
c \\
\displaystyle c-c\frac{a K_0  k}{2} \end{array} \right). \nonumber
\end{equation}
Choose $c=1$.

A general solution of (\ref{sys:PLLSysSawtooth}) has the form
\begin{align}
&\begin{cases}
\theta_\Delta(t) = e^{\left(\displaystyle -\frac{a K_0  k}{2} t\right)} \hskip0.1cm \left(c_1  + c_2\left(t + 1\right)\right), \\
y(t) = e^{\left(\displaystyle -\frac{a K_0  k}{2} t\right)} \hskip0.1cm \left(\displaystyle -\frac{a K_0  k}{2} c_1 + c_2\left(\displaystyle -\frac{a K_0  k}{2} t + 1 \displaystyle -\frac{a K_0  k}{2}\right)\right).
\end{cases} \label{eq:CommSolDegNode}
\end{align}
Similarly to \ref{subsec:StableNode} and \ref{subsec:StableFocus}, let us find the coefficients $c_1$ and $c_2$ for the solution of the Cauchy problem with the initial data $\theta_\Delta(0) = \frac{1}{k}$ and $y(0) = \frac{\sqrt{(aK_0 )^2 + 4{b}K_0 ({\pi}-\frac{1}{k})} - aK_0 }{2}$.
In this case we have:
\begin{align*}
&\begin{cases}
\displaystyle \frac{1}{k} = c_1 + c_2, \\
\displaystyle \frac{\sqrt{(aK_0 )^2 + 4{b}K_0 ({\pi}-\frac{1}{k})} - aK_0 }{2} = \displaystyle -\frac{a K_0  k}{2} c_1 + c_2\left(1 \displaystyle -\frac{a K_0  k}{2}\right),
\end{cases}
\end{align*}
\begin{align*}
&\begin{cases}
c_1 = \displaystyle \frac{1}{k} - c_2, \\
\displaystyle \frac{\sqrt{(aK_0 )^2 + 4{b}K_0 ({\pi}-\frac{1}{k})} - aK_0 }{2} = \displaystyle -\frac{a K_0  k}{2} \left(\frac{1}{k} - c_2\right) + c_2 -c_2\displaystyle \frac{a K_0  k}{2},
\end{cases}
\end{align*}
\begin{align*}
&\begin{cases}
c_1 = \displaystyle \frac{1}{k} - c_2, \\
\displaystyle \frac{\sqrt{(aK_0 )^2 + 4{b}K_0 ({\pi}-\frac{1}{k})}}{2} = c_2,
\end{cases}
\end{align*}
\begin{align*}
&\begin{cases}
c_1 = \displaystyle \frac{1}{k} - \displaystyle \frac{\sqrt{(aK_0 )^2 + 4{b}K_0 ({\pi}-\frac{1}{k})}}{2}, \\
c_2 = \displaystyle \frac{\sqrt{(aK_0 )^2 + 4{b}K_0 ({\pi}-\frac{1}{k})}}{2}.
\end{cases}
\end{align*}
Find $t_0$ such that $\theta_\Delta(t_0) = 0$:
\begin{align*}
&\begin{cases}
0 = e^{\left(\displaystyle -\frac{a K_0  k}{2} t_0\right)} \hskip0.1cm \left(c_1  + c_2\left(t_0 + 1\right)\right), \\
y(t_0) = e^{\left(\displaystyle -\frac{a K_0  k}{2} t_0\right)} \hskip0.1cm \left(\displaystyle -\frac{a K_0  k}{2} c_1 + c_2\left(\displaystyle -\frac{a K_0  k}{2} t_0 + 1 \displaystyle -\frac{a K_0  k}{2}\right)\right),
\end{cases}
\end{align*}
\begin{align*}
&\begin{cases}
0 = e^{\left(\displaystyle -\frac{a K_0  k}{2} t_0\right)} \hskip0.1cm \left(\left(c_1  + c_2 \right) + c_2 t_0 \right), \\
y(t_0) = e^{\left(\displaystyle -\frac{a K_0  k}{2} t_0\right)} \hskip0.1cm \left(\displaystyle -\frac{a K_0  k}{2} \left(c_1 + c_2\right) + c_2 \left(\displaystyle -\frac{a K_0  k}{2} t_0 + 1 \right)\right),
\end{cases}
\end{align*}
\begin{align*}
&\begin{cases}
0 = \displaystyle \frac{1}{k} + c_2 t_0, \\
y(t_0) = e^{\left(\displaystyle -\frac{a K_0  k}{2} t_0\right)} \hskip0.1cm \left(\displaystyle -\frac{a K_0 }{2} + c_2 \left(\displaystyle -\frac{a K_0  k}{2} t_0 + 1 \right)\right),
\end{cases}
\end{align*}
\begin{align*}
&\begin{cases}
t_0 = -\displaystyle \frac{1}{c_2 k}, \\
y(t_0) = e^{\left(\displaystyle -\frac{a K_0  k}{2} t_0\right)} \hskip0.1cm \left(\displaystyle -\frac{a K_0 }{2}  \displaystyle -\frac{a K_0  k c_2}{2}t_0 + c_2\right).
\end{cases}
\end{align*}
Finally, the expression for $S^\prime(\theta^s_{\rm eq})$ is as follows
\begin{equation}
S^\prime(\theta^s_{\rm eq}) = y(t_0) = c_2 \hskip0.1cm e^{\left(\displaystyle \frac{a K_0 }{2 c_2}\right)}, \label{rel:S0DegNode}
\end{equation}
where
\begin{equation}
c_2 = \displaystyle \frac{\sqrt{(aK_0 )^2 + 4{b}K_0 ({\pi}-\frac{1}{k})}}{2}. \nonumber
\end{equation}

\section*{\uppercase{Acknowledgments}}
 This work was supported by the Russian Scientific Foundation
 and Saint-Petersburg State University.
 The authors would like to thank Roland~E.~Best,
 the founder of the Best Engineering Company, Oberwil, Switzerland
 and the author of the bestseller on PLL-based circuits \cite{Best-2007}
 for valuable discussion.

%\section*{\uppercase{References}}
%\bibliographystyle{plainnat}
% \bibliography{bib_my,E:/Dropbox/bib/bib_pll,E:/Dropbox/bib/bib_nk,E:/Dropbox/bib/bib_full,E:/Dropbox/bib/bib_leonov,E:/Dropbox/bib/bib-2004-gly}
%\bibliography{bib_my,D:/Dropbox/bib/bib_pll,D:/Dropbox/bib/bib_nk,D:/Dropbox/bib/bib_full,D:/Dropbox/bib/bib_leonov,D:/Dropbox/bib/bib-2004-gly}
%\bibliography{C:/Dropbox/bib/bib_pll,C:/Dropbox/bib/bib_nk,C:/Dropbox/bib/bib_full,C:/Dropbox/bib/bib_leonov,C:/Dropbox/bib/bib-2004-gly}

\newcommand{\noopsort}[1]{} \newcommand{\printfirst}[2]{#1}
  \newcommand{\singleletter}[1]{#1} \newcommand{\switchargs}[2]{#2#1}

\end{document}